\documentclass[amsppt,11pt]{amsart}
\usepackage{amsmath}
\usepackage{amsfonts}
\usepackage{amssymb}
\usepackage[all]{xy}

\xyoption{dvips}

\newcommand{\ep}{{\varepsilon}}
\newcommand{\epi}{{\varepsilon_{\iota}}}
\newcommand{\intg}{{\mathbb Z}}
\newcommand{\cplx}{{\mathbb C}}

\newcommand{\fifi}{{\mathbb F}_q}
\newcommand{\fifisq}{{\mathbb F}_{q^2}}
\newcommand{\fifim}{{\mathbb F}_{q^m}}

\newcommand{\mltM}{{\mathbb M}}

\newcommand{\fldK}{{\mathbb K}}
\newcommand{\prt}{\mathcal{P}}
\newcommand{\bmu}{\boldsymbol{\mu}}
\newcommand{\blam}{\boldsymbol{\lambda}}
\newcommand{\bnu}{\boldsymbol{\nu}}
\newcommand{\bgamma}{\boldsymbol{\gamma}}
\newcommand{\balpha}{\boldsymbol{\alpha}}
\newcommand{\bbeta}{\boldsymbol{\beta}}

\newcommand{\spanning}{\text{-span}}

\newcommand{\ch}{\mathrm{ch}}
\newcommand{\hgt}{\mathrm{ht}}
\newcommand{\dd}{\displaystyle}
\newcommand{\Id}{\mathrm{Id}}
\newcommand{\One}{{\bf 1}}
\newcommand{\sss}{\scriptscriptstyle}
\newcommand{\bh}{\mathbf{h}}
\newcommand{\Hom}{\mathrm{Hom}}
\newcommand{\Ind}{\mathrm{Ind}}

\newcommand{\OneOne}{{1\hspace{-.14cm} 1}}
\newcommand{\oneone}{{1\hspace{-.11cm} 1}}

\newtheorem{thm}{Theorem}[section]
\newtheorem{lemma}{Lemma}[section]
\newtheorem{prop}{Proposition}[section]
\newtheorem{cor}{Corollary}[section]

\def\adots{\mathinner{\mkern2mu\raise0pt\hbox{.}  
\mkern2mu\raise4pt\hbox{.}\mkern1mu
\raise7pt\vbox{\kern7pt\hbox{.}}\mkern1mu}}

\numberwithin{equation}{section}

\allowdisplaybreaks[1]

\begin{document}

\bibliographystyle{amsplain}

\title{On the characteristic map of finite unitary groups}
\author{Nathaniel Thiem}
  \address{Department of Mathematics\\
           Stanford University\\
           Stanford, CA  94305-2125}
   \email{thiem@math.stanford.edu}
\author{C. Ryan Vinroot}
  \address{Department of Mathematics\\
           The University of Arizona\\
           617 N Santa Rita Ave\\
           Tucson, AZ  85721-0089}
   \email{vinroot@math.arizona.edu}

\maketitle

\section{Introduction}
In his seminal work \cite{green}, Green described a remarkable connection between the class functions of the finite general linear group ${\rm GL}(n, \fifi)$ and a generalization of the ring of symmetric functions of the symmetric group $S_n$.   In particular, Green defines a map, called the characteristic map, that takes irreducible characters to Schur-like symmetric functions, and recovers the character table of ${\rm GL}(n, \fifi)$ as the transition matrix between these Schur functions and Hall-Littlewood polynomials \cite[Chapter IV]{mac}.   Thus, we can use the combinatorics of the symmetric group $S_n$ to understand the representation theory of ${\rm GL}(n, \fifi)$.  Some of the implications of this approach include an indexing of irreducible characters and conjugacy classes of ${\rm GL}(n, \fifi)$ by 
multi-partitions and a formula for the degrees of the irreducible characters in terms of these partitions.

This paper describes the parallel story for the finite unitary group ${\rm U}(n, \fifisq)$ by collecting known results for this group and examining some applications of the unitary characteristic map.  Inspired by Green, Ennola \cite{ennolaconj,ennola} used results of Wall \cite{wall} to construct the appropriate ring of symmetric functions and characteristic map.  Ennola was able to prove that the analogous Schur-like functions correspond to an orthonormal basis for the class functions, and conjectured that they corresponded to the irreducible characters.   He theorized that the representation theory of ${\rm U}(n, \fifisq)$ should be deduced from the representation theory of ${\rm GL}(n, \fifi)$ by substituting ``$-q$'' for every occurrence of ``$q$''.  The general phenomenon of obtaining a polynomial invariant in $q$ for ${\rm U}(n, \fifisq)$ by this substitution has come to be known as ``Ennola duality''.  

Roughly a decade after Ennola made his conjecture, Deligne and Lusztig \cite{dellusz} constructed a family of virtual characters, called Deligne-Lusztig characters, to study the representation theory of arbitrary finite reductive groups.   Lusztig and Srinivasan \cite{luszsrin} then  computed an explicit decomposition of the irreducible characters of ${\rm U}(n, \fifisq)$ in terms of Deligne-Lusztig characters.  Kawanaka \cite{kawan2} used this composition to demonstrate that Ennola duality applies to Green functions, thereby improving results of Hotta and Springer \cite{hotspr} and finally proving Ennola's conjecture. 

This paper begins by describing some of the combinatorics and group theory associated with the finite unitary groups.  Section \ref{Preliminaries} defines the finite unitary groups, outlines the combinatorics of multi-partitions, and gives a description of some of the key subgroups.  Section \ref{Conjugacy} analyzes the conjugacy classes of ${\rm U}(n,\fifisq)$ and the Jordan decomposition of these conjugacy classes.

Section \ref{mapch} outlines the statement and development of the Ennola conjecture from two perspectives.  Both points of view define a map from a ring of symmetric functions to the character ring $C$ of ${\rm U}(n,\fifisq)$.   However, the first uses the multiplication for $C$ as defined by Ennola, and the second uses Deligne-Lusztig induction as the multiplicative structure of $C$.   The main results  are  

\vspace{.15cm}

\noindent \textbf{I.} (Theorem \ref{chDLcharacters}) The Deligne-Lusztig characters correspond to power-sum symmetric functions via the characteristic map of Ennola.

\vspace{.15cm}

\noindent \textbf{II.} (Corollary \ref{sameprod}) The multiplicative structure that Ennola defined on $C$ is Deligne-Lusztig induction.

\vspace{.15cm}

Section \ref{SCharacterDegrees} computes the degrees of the irreducible characters, and uses this result to evaluate various sums of character degrees   (see \cite[IV.6, Example 5]{mac} for the  ${\rm GL}(n, \fifi)$ analogue of this method). The main results are  

\vspace{.15cm}

\noindent\textbf{III.} (Theorem \ref{CharacterDegrees}) An irreducible $\chi^{\blam}$ character of ${\rm U}(m, \fifisq)$ corresponds to 
$$(-1)^{\lfloor m/2\rfloor+n(\blam)} s_{\blam}\quad \text{and}\quad \chi^{\blam}(1)=q^{n(\blam^\prime)} \frac{\dd\prod_{1\leq i\leq m}(q^i-(-1)^i)}{\dd\prod_{\Box\in \blam}(q^{\bh(\Box)} -(-1)^{\bh(\Box)})},$$
where $s_{\blam}$ is a Schur-like function, and both $n(\blam)$ and $\bh(\Box)$ are combinatorial statistics on the multi-partition $\blam$. 

\vspace{.15cm}

\noindent\textbf{IV.} (Corollary \ref{unsym}) If $\prt^\Theta_n$ indexes the irreducible characters $\chi^{\blam}$ of ${\rm U}(n, \fifisq)$, then
$$\sum_{\blam\in \prt_n^\Theta} \chi^{\blam}(1)=|\{g\in {\rm U}(n, \fifisq)\ \mid\ g\text{ symmetric}\}|.$$

\vspace{.15cm}

\noindent\textbf{V.} (Theorem \ref{evendegsum}) We give a subset $\mathcal{X}\subseteq \prt_{2n}^\Theta$ such that
$$\sum_{\blam\in\mathcal{X}} \chi^{\blam}(1)=(q+1)q^2(q^3+1)\cdots q^{2n-2}(q^{2n-1}+1) = \frac{|{\rm U}(2n, \fifisq)|}{|{\rm Sp}(2n, \fifi)|}.$$

Section \ref{Model} uses results by Ohmori \cite{ohm} and Henderson \cite{hender}  to adapt a model for the general linear group, found by Klyachko \cite{kl} and Inglis and Saxl \cite{ingsax}, to the finite unitary group.  The main result is

\vspace{.15cm}

\noindent\textbf{VI.} (Theorem \ref{model})  Let $U_m={\rm U}(m, \fifisq)$, where $q$ is odd, and let $\Gamma_{m}$ be the Gelfand-Graev character of $U_m$,  $\One$ be the trivial character of the finite symplectic group $Sp_{2r} = {\rm Sp}(2r, \fifi)$, and $R_L^G$ be the Deligne-Lusztig induction functor.  Then 
$$\sum_{0 \leq 2r \leq m} R_{U_{m-2r}\oplus U_{2r}}^{U_{m}}\bigl( \Gamma_{m-2r}\otimes \Ind_{Sp_{2r}}^{U_{2r}}(\One)\bigr)=\sum_{\blam\in \prt_m^\Theta} \chi^{\blam}.$$
That is, in the theorem of Klyachko, one may replace parabolic induction by Deligne-Lusztig induction to obtain a theorem for the unitary group. 

These results give considerable combinatorial control over the representation theory of the finite unitary group, and there are certainly more applications to these results than what we present in this paper.  Furthermore, this characteristic map gives some insight as to how a characteristic map might look in general type, using the invariant rings of other Weyl groups.

%
%
\section{Preliminaries}\label{Preliminaries}

%
%

\subsection{The unitary group and its underlying field $k$}\label{MandChars}

Let $\mathbb{K}=\bar{\mathbb{F}}_q$ the the algebraic closure of a finite field with $q$ elements and let $\mathbb{K}_m=\fifim$ denote the finite subfield with $q^m$ elements.   Let ${\rm GL}(n,\fldK)$ denote the general linear group over $\fldK$, and define Frobenius maps
\begin{equation}\begin{array}{r@{}c@{\ }c@{\ }c} F: & {\rm GL}(n,\fldK)  & \longrightarrow & {\rm GL}(n,\fldK) \\ & (a_{ij}) & \mapsto & (a_{ji}^q)^{-1},\end{array} \text{ and }  \begin{array}{r@{}c@{\ }c@{\ }c} F^\prime: & {\rm GL}(n,\fldK)  & \longrightarrow & {\rm GL}(n,\fldK)\\ & (a_{ij}) & \mapsto & (a_{n-j,n-i}^q)^{-1}.\end{array}
\end{equation}
Then the unitary group $U_n={\rm U}(n,\fldK_2)$ is given by
\begin{align} U_n&={\rm GL}(n,\fldK)^F=\{a\in {\rm GL}(n,\fldK)\ \mid\  F(a)=a\} \label{UnitaryDefinition}\\
			     &\cong {\rm GL}(n,\fldK)^{F^\prime} =\{a\in {\rm GL}(n,\fldK)\ \mid\ F^\prime(a)=a\}.\label{UnitaryDefinition2}
\end{align}
In fact, it follows from the Lang-Steinberg theorem (see, for example, \cite{Ca85}) that there exists $y\in {\rm GL}(n,\fldK_2)$ such that ${\rm GL}(n,\fldK)^F=y{\rm GL}(n,\fldK)^{F^\prime}y^{-1}$.

We define the multiplicative groups $\mltM_m$ as  
$$\mathbb{M}_m={\rm GL}(1,\fldK)^{F^m}=\{x\in \fldK \ \mid\ x^{q^m-(-1)^m}=1\}.$$  
Note that $\mltM_m \cong \fldK_m^{\times}$ only if $m$ is even.   We identify $\fldK^\times$ with the inverse limit ${\dd\lim_{\leftarrow} \mltM_m}$ with respect to the norm maps
$$ \begin{array}{rccc} N_{mr}:& \mltM_m & \longrightarrow & \mltM_r\\ & x & \mapsto & x x^{-q} \cdots x^{(-q)^{m/r-1}}\end{array},\quad \text{where $m,r\in \intg_{\geq 1}$ with $r\mid m$}.$$
If $\mltM_m^*$ is the group of characters of $\mltM_m$, then the direct limit  
$$\fldK^*=\lim_{\rightarrow} \mltM_m^*$$
gives the group of characters of $\fldK^\times$.  Let
$$\Theta=\{F\text{-orbits of $\fldK^{*}$}\}.$$ 

A polynomial $f(t)\in \fldK_2[t]$ is \emph{$F$-irreducible} if there exists an $F$-orbit $\{x_, x^{-q},\ldots, x^{(-q)^{d}}\}$ of $\fldK^\times$ such that
$$f(t)=(t-x)(t-x^{-q})\cdots (t-x^{(-q)^d}).$$
Let
\begin{equation}
\Phi=\{f\in \fldK_2[t]\ \mid\ f\text{ is $F$-irreducible}\}\overset{1-1}{\longleftrightarrow} \{F\text{-orbits of $\fldK^\times$}\}.
\end{equation}

The set $\Phi$ has an alternate description, as given in \cite{ennolaconj}.  For $f=t^d+a_1 t^{d-1} +\cdots + a_d\in \fldK_2[t]$ with $a_d\neq 0$, let 
$$ \tilde{f}=t^d + a_d^{-q}(a_{d-1}^q t^{d-1} +\cdots + a_1^q t+1),$$ 
which has the effect of applying $F$ to the roots of $f$ in $\fldK$.
Then a polynomial $f\in \fldK_2[t]$ is $F$-irreducible if and only if either
\begin{enumerate}
\item[(a)] $f$ is irreducible in $\fldK_2[t]$ and $\tilde{f}=f$ ($d$ must be odd in this case),  or
\item[(b)] $f(t) = h(t) \tilde{h}(t)$ where $h$ is irreducible in $\fldK_2[t]$ and $\tilde{h}(t) \neq h(t)$. 
\end{enumerate}

%
%

\subsection{Combinatorics of $\Phi$-partitions and $\Theta$-partitions}\label{Combinatorics}

Fix an ordering of $\Phi$ and $\Theta$, and let
$$\prt=\{\text{partitions}\} \quad \text{and}\quad \prt_n=\{\nu\in \prt\ \mid\ |\nu|=n\}.$$
Let $\mathcal{X}$ be either $\Phi$ or $\Theta$.    An \emph{$\mathcal{X}$-partition} $\bnu=(\bnu(x_1), \bnu(x_2),\ldots)$ is a sequence of partitions indexed by $\mathcal{X}$.   The \emph{size} of an $\mathcal{X}$-partition $\bnu$ is
\begin{equation}\label{part} ||\bnu||=\sum_{x\in \mathcal{X}} |x| |\bnu(x)|,\quad \text{where}\quad |x|=\left\{\begin{array}{ll} |x| & \text{if $\mathcal{X}= \Theta$,}\\ d(x) &  \text{if $\mathcal{X}= \Phi$,}\end{array}\right.
\end{equation}
$|x|$ is the size of the orbit $x\in \Theta$, and $d(x)$ is the degree of the polynomial $x\in \Phi$.  
Let 
\begin{equation}\label{Xpartitions}\prt^{\mathcal{X}}_n=\{\mathcal{X}\text{-partitions } \bnu\ \mid\ ||\bnu||=n\} ,\quad \text{and}\quad \prt^{\mathcal{X}}=\bigcup_{n=1}^{\infty} \prt^{\mathcal{X}}_n .\end{equation}

For $\bnu\in \prt^{\mathcal{X}}$, let
\begin{equation} \label{sigma}
n(\bnu)=\sum_{x\in \mathcal{X}} |x| n(\bnu(x)),\quad \text{where}\quad n(\nu)=\sum_{i=1}^{\ell(\nu)} (i-1)\nu_i.
\end{equation}
The \emph{conjugate} $\bnu^\prime$ of $\bnu$ is the $\mathcal{X}$-partition $\bnu^\prime=(\bnu(x_1)^\prime, \bnu(x_2)^\prime,\ldots )$, where $\nu^\prime$ is the usual conjugate partition for $\nu\in \prt$. 

The \emph{semisimple part} $\bnu_s$ of $\bnu=(\bnu(x_1),\bnu(x_2),\ldots)\in \prt_n^{\mathcal{X}}$ is
\begin{equation} \label{muSemisimple}\bnu_s=((1^{|\bnu(x_1)|}), (1^{|\bnu(x_2)|}),\ldots)\in \prt_n^{\mathcal{X}},\end{equation}
and the \emph{unipotent part} $\bnu_u$ of $\bnu\in \prt_n^{\mathcal{X}}$ is given by 
\begin{equation}\label{muUnipotent}\bnu_u(\OneOne) \quad \text{has parts} \quad \{|x|\bnu(x)_i \ \mid \ x \in \mathcal{X}, i = 1, \ldots, \ell(\bnu(x)) \}
\end{equation}
where
$$\OneOne=\left\{\begin{array}{ll} \{\One\} & \text{if $\mathcal{X}=\Theta$,}\\ t-1 & \text{if $\mathcal{X}=\Phi$,} \end{array}\right.$$
$\One$ is the trivial character in $\fldK^*$, and $\bnu_u(x) = \emptyset$ for $x \neq \OneOne$. 

\vspace{.15cm}

\noindent\textbf{Example.}  If $\bmu\in \prt^{\Phi}$ is given by
$$\bmu=\left(\xy<0cm,.3cm>\xymatrix@R=.3cm@C=.3cm{
	*={} & *={} \ar @{-} [l] & *={} \ar @{-} [l]  \\
	*={} \ar @{-} [u] &   *={} \ar @{-} [l] \ar@{-} [u] &   *={} \ar @{-} [l] \ar@{-} [u]   \\
*={} \ar @{-} [u] &   *={} \ar @{-} [l] \ar@{-} [u]  & *={} \ar @{-} [l] \ar@{-} [u]  }\endxy^{(f)}, \  \xy<0cm,.3cm>\xymatrix@R=.3cm@C=.3cm{
	*={} & *={} \ar @{-} [l] & *={} \ar @{-} [l] \\
	*={} \ar @{-} [u] &   *={} \ar @{-} [l] \ar@{-} [u]  &  *={} \ar @{-} [l] \ar@{-} [u] }\endxy^{(g)},\ 
\xy<0cm,.3cm>\xymatrix@R=.3cm@C=.3cm{
	*={} & *={} \ar @{-} [l] & *={} \ar @{-} [l] & *={} \ar @{-} [l] & *={} \ar @{-} [l]  \\
	*={} \ar @{-} [u] &   *={} \ar @{-} [l] \ar@{-} [u] &   *={} \ar @{-} [l] \ar@{-} [u]  &   *={} \ar @{-} [l] \ar@{-} [u] &   *={} \ar @{-} [l] \ar@{-} [u] \\ *={} \ar @{-} [u] &   *={} \ar @{-} [l] \ar@{-} [u] }\endxy^{(h)}\ \right),\quad\text{where $d(f)=1$, $d(g)=d(h)=2$,}$$
then $||\bmu||=1\cdot 4+2\cdot 2+2\cdot 5=18$, $n(\bmu)=1\cdot 2+ 2\cdot 0+ 2\cdot 1=4$,
$$\bmu_s=\left(\xy<0cm,.6cm>\xymatrix@R=.3cm@C=.3cm{
	*={} & *={} \ar @{-} [l] \\
	*={} \ar @{-} [u] &   *={} \ar @{-} [l] \ar@{-} [u]   \\
*={} \ar @{-} [u] &   *={} \ar @{-} [l] \ar@{-} [u]    \\
*={} \ar @{-} [u] &   *={} \ar @{-} [l] \ar@{-} [u]  \\
*={} \ar @{-} [u] &   *={} \ar @{-} [l] \ar@{-} [u]  }\endxy^{(f)}, \  
\xy<0cm,.3cm>\xymatrix@R=.3cm@C=.3cm{
	*={} & *={} \ar @{-} [l] \\
	*={} \ar @{-} [u] &   *={} \ar @{-} [l] \ar@{-} [u]  \\
*={} \ar @{-} [u] &   *={} \ar @{-} [l] \ar@{-} [u]   }\endxy^{(g)},\ 
\xy<0cm,.7cm>\xymatrix@R=.3cm@C=.3cm{
	*={} & *={} \ar @{-} [l]  \\
	*={} \ar @{-} [u] &   *={} \ar @{-} [l] \ar@{-} [u] \\ 
	*={} \ar @{-} [u] &   *={} \ar @{-} [l] \ar@{-} [u]  \\
*={} \ar @{-} [u] &   *={} \ar @{-} [l] \ar@{-} [u]  \\
*={} \ar @{-} [u] &   *={} \ar @{-} [l] \ar@{-} [u]  \\
*={} \ar @{-} [u] &   *={} \ar @{-} [l] \ar@{-} [u] }\endxy^{(h)}\ \right)\quad \text{and}\quad \bmu_u= \left(\xy<0cm,.7cm>\xymatrix@R=.3cm@C=.3cm{
	*={} & *={} \ar @{-} [l] & *={} \ar @{-} [l] & *={} \ar @{-} [l] & *={} \ar @{-} [l] & *={} \ar @{-} [l] & *={} \ar @{-} [l] & *={} \ar @{-} [l] & *={} \ar @{-} [l] \\
	*={} \ar @{-} [u] &   *={} \ar @{-} [l] \ar@{-} [u] &   *={} \ar @{-} [l] \ar@{-} [u]  &   *={} \ar @{-} [l] \ar@{-} [u] &   *={} \ar @{-} [l] \ar@{-} [u] &   *={} \ar @{-} [l] \ar@{-} [u] &   *={} \ar @{-} [l] \ar@{-} [u] &   *={} \ar @{-} [l] \ar@{-} [u] &   *={} \ar @{-} [l] \ar@{-} [u]  \\
	*={} \ar @{-} [u] &   *={} \ar @{-} [l] \ar@{-} [u]  & *={} \ar @{-} [l] \ar@{-} [u] & *={} \ar @{-} [l] \ar@{-} [u] & *={} \ar @{-} [l] \ar@{-} [u] \\
	*={} \ar @{-} [u] &   *={} \ar @{-} [l] \ar@{-} [u]  & *={} \ar @{-} [l] \ar@{-} [u]\\
	*={} \ar @{-} [u] &   *={} \ar @{-} [l] \ar@{-} [u]  & *={} \ar @{-} [l] \ar@{-} [u]\\
	*={} \ar @{-} [u] &   *={} \ar @{-} [l] \ar@{-} [u]  & *={} \ar @{-} [l] \ar@{-} [u] }\endxy^{(t-1)}\right).$$

\subsection{Levi subgroups and maximal tori}\label{LeviTorus}

Let $\mathcal{X}$ be either $\Phi$ or $\Theta$ as in Section \ref{Combinatorics}.  

For $\bnu\in \prt^{\mathcal{X}}_n$, let
\begin{equation}\label{Lmu} L_{\bnu}=\bigoplus_{x\in \mathcal{X}_{\bnu}} L_{\bnu}(x),\quad\text{where}\quad  \mathcal{X}_{\bnu}=\{x\in \mathcal{X}\ \mid\ \bnu(x)\neq \emptyset\},\end{equation}
and for $x\in \mathcal{X}_{\bnu}$,
\begin{equation} \label{Lmuf}
L_{\bnu}(x)=\left\{\begin{array}{ll} {\rm U}(|\bnu(x)|, \fldK_{2|x|})& \text{if $|x|$ is odd,}\\ {\rm GL}(|\bnu(x)|, \fldK_{|x|}) &\text{if $|x|$ is even.}\end{array}\right.
\end{equation}
Then $L_{\bnu}$ is a Levi subgroup of $U_n={\rm U}(n,\fldK_2)$ (though not uniquely determined by $\bnu$).  
The Weyl group 
\begin{equation}\label{LeviWeyl} W_{\bnu}=\bigoplus_{x\in \mathcal{X}_{\bnu}} S_{|\bnu(x)|},\end{equation}
 of $L_{\bnu}$ has conjugacy classes indexed by
\begin{equation} \label{SemisimpleMatch} \prt^{\bnu}_s=\{\bgamma\in \prt^{\mathcal{X}}\ \mid\ \bgamma_s=\bnu_s\},\end{equation}
and the size of the conjugacy class $c_{\bgamma}$ is
\begin{equation}\label{zbgamma}
|c_{\bgamma}|=\frac{|W_{\bgamma}|}{z_{\bgamma}},\quad \text{where}\quad z_{\bgamma}=\prod_{x\in \mathcal{X}} z_{\bgamma(x)}\quad\text{and}\quad z_\gamma=\prod_{i=1}^{\ell(\gamma)} i^{m_i}m_i!,
\end{equation}
for $\gamma=(1^{m_1}2^{m_2}\cdots)\in \prt$.

For every $\nu=(\nu_1,\nu_2,\ldots, \nu_\ell)\in \prt_n$ there exists a maximal torus (unique up to isomorphism) $T_\nu$ of $U_n$ such that 
$$T_\nu\cong \mltM_{\nu_1} \times \mltM_{\nu_2} \times\cdots\times \mltM_{\nu_\ell}.$$ 
For every $\bgamma\in \prt_s^{\bmu}$, there exists a maximal torus (unique up to isomorphism) $T_{\bgamma}\subseteq L_{\bnu}$ such that
\begin{equation}\label{Torus} T_{\bgamma}=\bigoplus_{x\in \mathcal{X}_{\bnu}} T_{\bgamma}(x),\quad\text{where} \quad T_{\bgamma}(x)\cong \mltM_{|x|\bgamma(x)_1}\times\cdots \times \mltM_{|x|\bgamma(x)_{\ell}}.
\end{equation}
Note that as a maximal torus of $U_n$, the torus $T_{\bgamma}\cong T_{\bgamma_u(\oneone)}$.

%
%

\section{Conjugacy classes and Jordan decomposition}\label{Conjugacy}

Let $U_n = {\rm U}(n, \fldK_2)$ as in (\ref{UnitaryDefinition}).

\begin{prop}
The conjugacy classes $c_{\bmu}$ of $U_n$ are indexed by $\bmu\in \prt_n^\Phi$.
\end{prop}
\begin{proof} 
For any $g \in U_n$, the monic minimal polynomial $f \in \fldK_2[t]$ of $g$ satisfies $f(t) = \tilde{f}(t)$, and in fact every invariant factor of $g$ satisfies this.  Every polynomial of this kind may be factored into $F$-irreducibles.   Wall \cite{wall} proved that an element $g \in {\rm GL}(n, \fldK_2)$ is an element of $U_n$ if and only if for every irreducible $f$, the polynomials $f$ and $\tilde{f}$ have the same multiplicities as elementary divisors of $g$.  

Each $g \in U_n$ acts on the vector space $\fldK_2^n$ and thus gives it a $\fldK_2[t]$-module structure, where $t \cdot v = gv$ for each $v \in \fldK_2^n$.  It is a result of Wall \cite{wall} that two elements of $U_n$ are conjugate in $U_n$ if they are conjugate in $G_n = {\rm GL}(n, \fldK_2)$.  So, two elements of $U_n$ are conjugate if and only if they induce the same $\fldK_2[t]$-module structure on $\fldK_2^n$.  For a conjugacy class $c$ of $U_n$, we call this $\fldK_2[t]$-module $V_c$.  Since the minimal polynomial $f$ of $g \in U_n$ satisfies $f(t) = \tilde{f}(t)$ and annihilates $V_c$, we have that the conjugacy classes of $U_n$ are parametrized by $\fldK_2[t]$-modules whose annihilators are generated by such polynomials.

It follows that $V_c$ is the direct sum of cyclic modules of the form $\fldK_2[t]/(h)^m$, where $h$ is irreducible, and $(\tilde{h})^m$ appears in the direct sum whenever $(h)^m$ occurs.  That is, when $h^m$ and $\tilde{h}^m$ are distinct, we may pair them together.  In this way, to each $f \in \Phi$ we have assigned a partition $\bmu(f) = (\bmu(f)_1, \bmu(f)_2, \ldots)$ such that
$$ V_c \cong \bigoplus_{f, i} \fldK_2[t]/(f)^{\bmu(f)_i} .$$
Since $f \in \Phi$, then either $f$ is irreducible, or $f = g \tilde{g}$ where $g$ is irreducible.  So $V_c$ determines a $\Phi$-partition $\bmu$  of size $n$.

Conversely, every $\Phi$-partition $\bmu$ of size $n$ determines a conjugacy class of $U_n$.  
\end{proof} 

For $r\in \intg_{\geq 0}$, let
\begin{equation}\label{psiProduct} \psi_r(x)=\prod_{i=1}^r(1-x^i).\end{equation}
\begin{prop}[Wall] \label{wall}  Let $g\in c_{\bmu}$.  The order $a_{\bmu}$ of the centralizer $g$ in $U_n$ is 
$$a_{\bmu}=(-1)^{||\bmu||}\prod_{f\in \Phi}a_{\bmu(f)} \bigl((-q)^{d(f)}\bigr),\text{ where }
a_{\mu}(x) = x^{|\mu| + 2n(\mu)} \prod_j \psi_{m_j} (x^{-1}),$$
for $\mu=(1^{m_1}2^{m_2}3^{m_3}\cdots)\in \prt$. 
\end{prop}

\begin{proof}
It follows from the results of Wall \cite[Section 2.6]{wall}, along with the description of $F$-irreducible polynomials given in Section 2 above, that
$$ a_{\bmu}=\prod_{f \in \Phi}q^{d(2\sum_{i<j} im_i m_j + \sum_i (i-1) m_i^2)} \prod_j \Big( q^{dm_j (m_j - 1)/2} \prod_{i = 1}^{m_j} (q^{di} - (-1)^{di}) \Big).$$
where $d=d(f)$ and $\bmu(f)=(1^{m_1}2^{m_2}\cdots)$.
From
$$ \prod_{i = 1}^{m_j} (q^{di} - (-1)^i) = q^{dm_j (m_j + 1)/2} \prod_{i = 1}^{m_j} (1 - (1/(-q)^d)^i), $$
we obtain
$$a_{\bmu}=\prod_{f\in \Phi} q^{d(2\sum_{i<j} im_i m_j + \sum_i i m_i^2)} \prod_j \psi_{m_j} \big( (-q)^{-d} \big). $$
Note that for any partition $\mu=(1^{m_1}2^{m_2}\cdots)\in \prt$, it follows from
$$ \sum_{j \leq i+1} m_j = \mu_{i+1}^{\prime}\qquad \text{and}\qquad m_i = \mu_i^{\prime} - \mu_{i+1}^{\prime},$$
that
$$2\sum_{i<j} im_i m_j + \sum_i i m_i^2= |\mu| + 2n(\mu).$$
Thus 
\begin{align*}
  a_{\bmu} & =  \prod_{f\in \Phi}(q^d)^{|\bmu(f)| + 2n(\bmu(f))} \prod_j \psi_{m_j} \big( (-q)^{-d} \big)\\
		& = \prod_{f\in \Phi}(-1)^{d|\bmu(f)|} ((-q)^d)^{|\bmu(f)| + 2n(\bmu(f))} \prod_j \psi_{m_j} \big( (-q)^{-d} \big)\\
		&=  (-1)^{||\bmu||}\prod_{f\in \Phi} a_{\bmu(f)}\bigl((-q)^{d(f)}\bigr),
		\end{align*}
as desired.
\end{proof}

For $\bmu\in \prt^\Phi$, let $L_{\bmu}$ be as in (\ref{Lmu}).  Note that $|L_{\bmu}|=a_{\bmu_s}.$

\begin{lemma} \label{JordanDecomposition} Suppose $g\in c_{\bmu}$ with Jordan decomposition $g=su$.  Then 
\begin{enumerate}
\item[(a)]  $s\in c_{\bmu_s}$ and $u\in c_{\bmu_u}$, where $\bmu_s$ and $\bmu_u$ are as in (\ref{muSemisimple}) and (\ref{muUnipotent}),
\item[(b)] the centralizer $C_{U_n}(s)$ of $s$ in $U_n$ is isomorphic to $ L_{\bmu}.$
\end{enumerate}
\end{lemma}
\begin{proof}
(a)  Suppose $f=t^d-a_{d-1}t^{d-1}-\cdots a_1t-a_0\in \Phi$ is irreducible (so $d$ is odd).  Define
\begin{equation} \label{CompanionMatrix} J(f)=\left(\begin{array}{ccccc} 0 & 1 & 0 & \cdots & 0\\ 0 & 0 & 1 & \ddots & \vdots\\ \vdots & & \ddots & \ddots & 0\\ 0 & 0 & \cdots & 0 & 1\\ a_0 & a_1 & a_2 & \cdots & a_{d-1}\end{array}\right)\in {\rm GL}(d, \fldK_2).\end{equation}
For $\nu=(\nu_1,\nu_2,\ldots,\nu_\ell)\in \prt$, let
$$J_\nu(f)=J_{\nu_1}(f)\oplus\cdots\oplus J_{\nu_\ell}(f), \text{ where } J_m(f)=\left(\begin{array}{cccc} J(f) & \Id_{d} & & 0\\ 0 & J(f) & \ddots & \\ \vdots & \ddots & \ddots & \Id_{d} \\ 0 & \cdots & 0 & J(f)\end{array}\right)$$
is an $md\times md$ matrix and $\Id_d$ is the $d\times d$ identity matrix.  Note that
$$J_m(f)=J_{(1^m)}(f)u_{(m)}(f),\quad  \text{where}\quad u_{(m)}(f)=\left(\begin{array}{cccc} \Id_d & J(f)^{-1} & & 0\\  & \Id_d & \ddots & \\ & & \ddots & J(f)^{-1} \\ 0 & & &\Id_d\end{array}\right)$$
is the Jordan decomposition of $J_m(f)$ in ${\rm GL}(d, \fldK_2)$.  

Similarly, if $f=h\tilde{h}\in \Phi$ (or $d(f)$ is even), then define
\begin{align*} J(f)&=J(h)\oplus J(\tilde{h})\\
J_m(f)&=J_m(h)\oplus J_m(\tilde{h})\\
J_\nu(f)&=J_\nu(h)\oplus J_\nu(\tilde{h}),\end{align*}
so
$$J_m(f)=\bigl(J_{(1^m)}(h)\oplus J_{(1^m)}(\tilde{h})\bigr)\bigl(u_{(m)}(h)\oplus u_{(m)}(\tilde{h})\bigr)$$
is the Jordan decomposition of $J_m(f)$.
If $g\in c_{\bmu}$, then $g$ is conjugate to
$$J_{\bmu}=\bigoplus_{f\in \Phi} J_{\bmu(f)}(f)$$
in ${\rm GL}(||\bmu||,\fldK_2)$.   Claim (a) follows from the uniqueness of the Jordan decomposition.

(b)   In ${\rm GL}(n, \fldK_2)$,
$$C_{{\rm GL}(n, \fldK_2)}(J_{(1^m)}(f))\cong \left\{\begin{array}{ll} {\rm GL}(m, \fldK_{2d(f)}) &\text{if $d(f)$ is odd,} \\ {\rm GL}(m, \fldK_{d(f)})\oplus {\rm GL}(m, \fldK_{d(f)}) & \text{if $d(f)$ is even.}\end{array}\right.$$
If $d(f)$ is odd and $s=xJ_{(1^m)}(f)x^{-1}$ for some $x\in {\rm GL}(n, \fldK_2)$, then 
$$(x{\rm GL}(m, \fldK_{2d(f)})x^{-1})^{F}= {\rm GL}(m, \fldK_{2d(f)})^{x\circ  F}\cong {\rm GL}(m, \fldK_{2d(f)})^{F}={\rm U}(m, \fldK_{2d(f)}).$$

Suppose $d(f)$ is even so that $f=h\tilde{h}$, and let $y\in {\rm GL}(n,\fldK_2)$ such that ${\rm GL}(n,\fldK)^{F}=y{\rm GL}(n,\fldK)^{F^\prime}y^{-1}$ (see the comment after (\ref{UnitaryDefinition})).  The element $J(f)$ is conjugate to
$$J(h)\oplus (F^\prime(J(h)))\in {\rm GL}(n, \fldK_2)^{F^\prime},$$
whose centralizer in ${\rm GL}(n, \fldK_2)^{F^\prime}$ is
$$\{g\oplus F^\prime(g)\ \mid\ g\in {\rm GL}(m, \fldK_{d(f)})\}\cong {\rm GL}(m, \fldK_{d(f)}).\qedhere$$
\end{proof}

%
%

\section{The Ennola Conjecture}\label{mapch}

%
%

\subsection{The characteristic map}\label{SymmetricFunctions}

Let $X=\{X_1,X_2,\ldots\}$ be an infinite set of variables and let $\Lambda(X)$ be the graded $\cplx$-algebra of symmetric functions in the variables $\{X_1,X_2,\ldots\}$.  For $\nu=(\nu_1,\nu_2,\ldots,\nu_\ell)\in \prt$, the \emph{power-sum symmetric function} $p_\nu(X)$ is
$$p_\nu(X)=p_{\nu_1}(X)p_{\nu_2}(X)\cdots p_{\nu_\ell}(X),\quad \text{where}\quad p_m(X)=X_1^m+X_2^m+\cdots.$$
The irreducible characters $\omega^\lambda$ of $S_n$ are indexed by $\lambda\in\prt_n$.  Let $\omega^\lambda(\nu)$ be the value of $\omega^\lambda$ on a permutation with cycle type $\nu$.  
For $\lambda\in \prt$, the \emph{Schur function} $s_\lambda(X)$ is given by
\begin{equation} \label{SchurToPower} s_\lambda(X)=\sum_{\nu\in \prt_{|\lambda|}}  \omega^\lambda(\nu)z_{\nu}^{-1} p_\nu(X),\quad\text{where}\quad  z_{\nu} = \prod_{i \geq 1} i^{m_i} m_i! \end{equation}
is the order of the centralizer in $S_n$ of the conjugacy class corresponding to $\nu=(1^{m_1}2^{m_2}\cdots)\in \prt$.
Let $t\in \cplx$.  For $\mu\in \prt$, the \emph{Hall-Littlewood symmetric function} $P_\mu(X;t)$ is given by
\begin{equation}\label{SchurToHall}
s_\lambda(X)=\sum_{\mu\in \prt_{|\lambda|}} K_{\lambda\mu}(t) P_\mu(X;t),
\end{equation}
where $K_{\lambda\mu}(t)$ is the Kostka-Foulkes polynomial (as in \cite[III.6]{mac}). For $\nu,\mu\in \prt_n$, the \emph{classical Green function} $Q_\nu^\mu(t)$ is given by
\begin{equation}\label{PowerToHall} p_\nu(X)=\sum_{\mu\in \prt_{|\nu|}} Q_\nu^\mu(t^{-1})t^{n(\mu)} P_\mu(X;t).\end{equation}
As a graded ring,
\begin{align*} \Lambda(X) &=\cplx\spanning\{p_\nu(X)\ \mid\ \nu\in \prt\}\\
					 &=\cplx\spanning\{s_\lambda(X)\ \mid\ \lambda\in \prt\}\\
					 &=\cplx\spanning\{P_\mu(X;t)\ \mid\ \mu\in \prt\}.
\end{align*}

For every $f \in \Phi$, fix a set of independent variables $X^{(f)}=\{X_1^{(f)},X_2^{(f)},\ldots\}$, and for any symmetric function $h$, we let $h(f) = h(X^{(f)})$ denote the symmetric function in the variables $X^{(f)}$.    
Let 
$$\Lambda=\cplx\spanning\{P_{\bmu}\ \mid\ \bmu\in \prt^\Phi\},\quad \text{where}\quad P_{\bmu} =(-q)^{-n(\bmu)} \prod_{f \in \Phi} P_{\bmu(f)} (f;(-q)^{d(f)}).$$
Then
$$\Lambda=\bigoplus_{n\geq 0}\Lambda_n,\quad \text{where}\quad \Lambda_n=\cplx\spanning\{P_{\bmu} \mid\ ||\bmu||=n\},$$
makes $\Lambda$ a graded $\cplx$-algebra.
Define a Hermitian inner product on $\Lambda$ by
$$ \langle P_{\bmu}, P_{\bnu} \rangle = a_{\bmu}^{-1} \delta_{\bmu\bnu}. $$

For each $\varphi\in \Theta$ let $Y^{(\varphi)}=\{Y_1^{(\varphi)},Y_2^{(\varphi)},\ldots\}$ be an infinite variable set, and for a symmetric function $h$, let $h(\varphi)= h(Y^{(\varphi)})$.   Relate symmetric functions in the $X$ variables to symmetric functions in the $Y$ variables via the transform
\begin{equation} \label{transform}
p_n(\varphi)=(-1)^{n|\varphi|-1}\sum_{x\in \mltM_{n|\varphi|}} \xi(x) p_{n|\varphi|/d(f_x)} (f_x),
\end{equation}
where $\varphi\in \Theta$, $\xi\in \varphi$, and $f_x\in \Phi$ satisfies $f_x(x)=0$. 

Then
\begin{equation}\Lambda=\cplx\spanning\{s_{\blam}\ \mid\ \blam\in \prt^\Theta\},\quad\text{where}\quad s_{\blam} = \prod_{\varphi \in \Theta} s_{\blam(\varphi)}(\varphi).\end{equation}

Let $C_n$ denote the set of complex-valued class functions of the group $U_n$, and for $||\bmu|| = n$, let $\pi_{\bmu}: U_n \rightarrow \cplx$ be given by
$$\pi_{\bmu}(u)=\left\{\begin{array}{ll} 1 & \text{if $u\in c_{\bmu}$},\\ 0 &\text{otherwise},\end{array}\right.\qquad \text{where $u\in U_n$}.$$ 
Then the $\pi_{\bmu}$ form a $\cplx$-basis for $C_n$.  By Proposition \ref{wall}, the usual inner product on class functions of finite groups, $\langle \cdot, \cdot \rangle:C_n\times C_n\rightarrow \cplx$, satisfies
$$\langle \pi_{\bmu},\pi_{\blam}\rangle= a_{\bmu}^{-1} \delta_{\bmu \blam}.$$

For $\alpha_i \in C_{n_i}$, Ennola \cite{ennola} defined a product $\alpha_1 \star \alpha_2 \in C_{n_1 + n_2}$, which takes the following value on the conjugacy class $c_{\blam}$:
$$ \alpha_1 \star \alpha_2 (c_{\blam}) =  \sum_{||\bmu_i|| = n_i} g_{\bmu_1 \bmu_2}^{\blam} \alpha_1 (c_{\bmu_1}) \alpha_2 (c_{\bmu_2}),$$
where $g_{\bmu_1 \bmu_2}^{\blam}$ is the product of Hall polynomials (see \cite[Chapter II]{mac})
$$ g_{\bmu_1 \bmu_2}^{\blam} = \prod_{f \in \Phi} g_{\bmu_1(f) \bmu_2(f)}^{\blam(f)} ((-q)^{d(f)}). $$
Extend the inner product to 
$$ C = \bigoplus_{n \geq 0} C_n, $$
by requiring the components $C_n$ and $C_m$ to be orthogonal for $n \neq m$.  This gives $C$ a graded $\cplx$-algebra structure.  The {\it characteristic map} is
$$\begin{array}{rccl} {\rm ch}:& C & \longrightarrow & \Lambda\\ &\pi_{\bmu} & \mapsto & P_{\bmu},\qquad \text{for $\bmu\in \prt^\Phi$}.\end{array}$$

\begin{prop} \label{charmap1} Let multiplication in the character ring $C$ of $U_n$ be given by $\star$.  Then the characteristic map $\ch: C \rightarrow \Lambda$ is an isometric isomorphism of graded $\cplx$-algebras.
\end{prop}
\begin{proof} It is clear that ${\rm ch}$ is a linear graded isometric isomorphism of vector spaces.  It is a ring homomorphism because 
\begin{align*}
\ch(\pi_{\bmu}\pi_{\bnu})&=\sum_{\blam\in \prt^\Phi} g_{\bmu\bnu}^{\blam}  \ch(\pi_{\blam})\\
&= \sum_{\blam \in \prt^\Phi} g_{\bmu\bnu}^{\blam} P_{\blam}\\
&= P_{\bmu}P_{\bnu} & & \text{(by \cite[III.3.6]{mac})}\\
&  = \ch(\pi_{\bmu})\ch(\pi_{\bnu}). & & \qedhere
\end{align*}
\end{proof}

%
%

\subsection{The  Ennola conjecture I}\label{EnnolaConj}

Following the work of Green \cite{green} on the general linear group, Ennola was able to obtain the following result.  We follow the proof in Macdonald \cite[IV.4]{mac} on the general linear group case.

\begin{prop}[Ennola] \label{orthog}
The set $\{ s_{\blam} \;\, | \;\, \blam \in \prt^{\Theta} \}$ is an orthonormal basis for $\Lambda$.
\end{prop}
\begin{proof}
For $u \in \Lambda$, if
$$u=\sum_{\bmu \in \prt^{\Phi}} u_{\bmu} P_{\bmu}\qquad \text{then let}\qquad \bar{u}=\sum_{\bmu\in \prt^{\Phi}} \bar{u}_{\bmu} P_{\bmu},$$
where $\bar{x}\in \cplx$ is the complex conjugate of $x\in \cplx$.  It follows from \cite[I.4.3]{mac} that
\begin{equation}
\log\left( \sum_{\blam\in \prt^\Theta} s_{\blam} \otimes \bar{s}_{\blam}\right) 
=\sum_{n\geq 1}\frac{1}{n}\sum_{\xi\in \mltM_n^*} p_{n/|\varphi_\xi|}(\varphi_\xi)\otimes \overline{p_{n/|\varphi_\xi|}(\varphi_\xi)}. \label{A1}
\end{equation} 
where $\varphi_\xi\in \Theta$ is the orbit such that $\xi\in \varphi_\xi$.
Apply the transform (\ref{transform}) and orthogonality of characters to obtain
\begin{equation}
\sum_{\xi\in \mltM_n^*} p_{n/|\varphi_\xi|}(\varphi_\xi)\otimes \overline{p_{n/|\varphi_\xi|}(\varphi_\xi)} = |\mltM_n|\sum_{x \in \mltM_n} p_{n/d(f_x)}(f_x) \otimes p_{n/d(f_x)}(f_x). \label{A2}
\end{equation}
From (\ref{A1}) and (\ref{A2}), we have
\begin{align}
\log\left( \sum_{\blam \in \prt^\Theta} s_{\blam} \otimes \bar{s}_{\blam}\right) & = \sum_{n \geq 1} \frac{(q^{n} - (-1)^{n})}{n} \sum_{x \in \mltM_n} p_{n/d(f_x)}(f_x) \otimes p_{n/d(f_x)}(f_x)\notag\\
 & = \sum_{n \geq 1} \frac{1}{n} \sum_{f \in \Phi} (q^{nd(f)} - (-1)^{nd(f)}) p_n(f) \otimes p_n(f). \label{A}
\end{align}

On the other hand, note that by \cite[III.4.4]{mac} the $P_\nu(t)$ satisfy a Cauchy kernel formula
$$\sum_{\nu\in\prt} t^{|\nu|+2n(\nu)}a_\nu(t^{-1}) P_\nu(x;t)P_\nu(y;t)=\prod_{i,j} \frac{1-tx_iy_j}{1-x_iy_j}.$$
With substitutions $t\mapsto (-q)^{-d(f)}$ and $x_i\mapsto q^{d(f)}x_i$ we obtain
\begin{align*}
\prod_{i,j} &\frac{1-(-1)^{d(f)}x_iy_j}{1-q^{d(f)}x_iy_j}\\
&=\sum_{\nu\in\prt} (-1)^{|\nu|d(f)}q^{-2n(\nu)d(f)}a_\nu((-q)^{d(f)}) P_\nu(x;(-q)^{-d(f)})P_\nu(y;(-q)^{-d(f)}).
\end{align*}
Make the further substitution $x_i\mapsto X_i^{(f)}\otimes 1$, $y_j\mapsto 1\otimes X_j^{(f)}$ and take the product on both sides over all $f\in \Phi$ to get
$$\prod_{f\in \Phi}\prod_{i,j} \frac{1-(-1)^{d(f)}X_i^{(f)}\otimes X_j^{(f)}}{1-q^{d(f)}(X_i^{(f)}\otimes X_j^{(f)})}=\sum_{\bmu \in \prt^{\Phi}} a_{\bmu} P_{\bmu}\otimes P_{\bmu}.$$
Take logarithms of both sides to obtain
\begin{align}
\log\biggl(&\sum_{\bmu \in \prt^{\Phi}} a_{\bmu} P_{\bmu}\otimes P_{\bmu} \biggr) \notag\\
 & = \sum_{f\in \Phi\atop i,j}  \biggl(\log(1-(-1)^{d(f)}X_i^{(f)}\otimes X_j^{(f)})-\log(1-q^{d(f)}X_i^{(f)}\otimes X_j^{(f)})\biggr)\notag\\
&=\sum_{n\geq 1} \frac{1}{n} \sum_{f\in \Phi} \sum_{i,j} (q^{d(f)n}-(-1)^{d(f)n}) \bigl(X_i^{(f)}\otimes X_j^{(f)}\bigr)^n\notag\\
&=\sum_{n\geq 1} \frac{1}{n} \sum_{f\in \Phi} (q^{nd(f)}-(-1)^{nd(f)}) p_n(f)\otimes p_n(f).\label{B}
\end{align}
We may conclude from (\ref{A}) and (\ref{B}) that 
\begin{equation} \sum_{\blam\in \prt^\Theta} s_{\blam} \otimes \bar{s}_{\blam}=
\sum_{\bmu \in \prt^{\Phi}} a_{\bmu} P_{\bmu}\otimes P_{\bmu}.\label{C}\end{equation}
Suppose
$$s_{\blam}=\sum_{\bmu\in \prt^\Phi} u_{\blam\bmu}P_{\bmu}=\sum_{\bmu\in \prt^\Phi} v_{\blam\bmu} a_{\bmu} P_{\bmu}.$$
Compare coefficients of $a_{\bmu} P_{\bmu}\otimes P_{\bmu}$ on both sides of (\ref{C}) to obtain
$$\sum_{\blam\in \prt^\Phi} u_{\blam\bmu} \bar{v}_{\blam\bnu}=\delta_{\bmu\bnu}.$$
It follows that $\langle s_{\blam},s_{\bnu}\rangle=\delta_{\blam\bnu}$, as desired.
\end{proof}

Now let $\chi^{\blam} \in R$ be class functions so that $\chi^{\blam}(1) > 0$ and $\ch(\chi^{\blam}) = \pm s_{\blam}$.  Ennola conjectured that $\{ \chi^{\blam} \;\, | \;\, \blam \in \prt^{\Theta}_n \}$ is the set of irreducible characters of $U_n$.  He pointed out that if one could show that the product $\star$ takes virtual characters to virtual characters, then the conjecture would follow.  There is no known direct proof of this fact, however.  Significant progress on Ennola's conjecture was only made after the work of Deligne and Lusztig \cite{dellusz}.

%
%
\subsection{Deligne-Lusztig Induction}\label{Deligne-Lusztig}

Let $T_\nu\cong \mltM_{\nu_1} \times \cdots \times \mltM_{\nu_\ell}$ be a maximal torus of $U_n$.   If $t\in T_\nu$, then $t$ is conjugate to 
$$J_{(1^{m_1})}(f_1)\oplus\cdots \oplus J_{(1^{m_\ell})}(f_\ell),\quad \text{where $f_i\in \Phi$, $m_id(f_i)=\nu_i$.}$$
Define $\bgamma_t\in \prt^\Phi$ by
\begin{equation}\label{gammat} \bgamma_t(f) \qquad\text{has parts}\qquad \{m_i\ \mid\  f_i=f\}.\end{equation}
Note that $(\bgamma_t)_u(t-1)=\nu$, but in general $t\notin c_{\bgamma_t}$.

\vspace{.25cm}

\noindent\textbf{Example.} If $t\in T_{(4,4,2,2,1)}$ and $t$ is conjugate to
$$J_{\xymatrix@R=.15cm@C=.15cm{*={} & *={}\ar @{-} [l] \\ *={} \ar @{-} [u] & *={}\ar @{-} [l] \ar @{-} [u]\\*={} \ar @{-} [u] & *={}\ar @{-} [l] \ar @{-} [u]\\ *={} \ar @{-} [u] & *={}\ar @{-} [l] \ar @{-} [u]\\ *={} \ar @{-} [u] & *={}\ar @{-} [l] \ar @{-} [u]} }(t-1)\oplus J_{\xymatrix@R=.15cm@C=.15cm{*={} & *={}\ar @{-} [l] \\ *={} \ar @{-} [u] & *={}\ar @{-} [l] \ar @{-} [u]\\*={} \ar @{-} [u] & *={}\ar @{-} [l] \ar @{-} [u]}}(t^2+1)\oplus J_{\xymatrix@R=.15cm@C=.15cm{*={} & *={}\ar @{-} [l] \\ *={} \ar @{-} [u] & *={}\ar @{-} [l] \ar @{-} [u]\\*={} \ar @{-} [u] & *={}\ar @{-} [l] \ar @{-} [u]}}(t-1) \oplus J_{\xymatrix@R=.15cm@C=.15cm{*={} & *={}\ar @{-} [l] \\ *={} \ar @{-} [u] & *={}\ar @{-} [l] \ar @{-} [u]}}(t^2+1)\oplus J_{\xymatrix@R=.15cm@C=.15cm{*={} & *={}\ar @{-} [l] \\ *={} \ar @{-} [u] & *={}\ar @{-} [l] \ar @{-} [u]}}(t-1),$$
then
$$\bgamma_t=\left(\xy<0cm,.5cm>\xymatrix@R=.3cm@C=.3cm{
	*={} & *={} \ar @{-} [l] & *={} \ar @{-} [l] & *={} \ar @{-} [l]  & *={} \ar @{-} [l] \\
	*={} \ar @{-} [u] &   *={} \ar @{-} [l] \ar@{-} [u] &   *={} \ar @{-} [l] \ar@{-} [u] &   *={} \ar @{-} [l] \ar@{-} [u]  &   *={} \ar @{-} [l] \ar@{-} [u]  \\
*={} \ar @{-} [u] &   *={} \ar @{-} [l] \ar@{-} [u] &   *={} \ar @{-} [l] \ar@{-} [u]  \\
*={} \ar @{-} [u] &   *={} \ar @{-} [l] \ar@{-} [u]   }\endxy^{(t-1)}, \  \xy<0cm,.3cm>\xymatrix@R=.3cm@C=.3cm{
	*={} & *={} \ar @{-} [l] & *={} \ar @{-} [l] \\
	*={} \ar @{-} [u] &   *={} \ar @{-} [l] \ar@{-} [u]  &  *={} \ar @{-} [l] \ar@{-} [u] \\
	*={} \ar @{-} [u] &   *={} \ar @{-} [l] \ar@{-} [u] 
	 }\endxy^{(t^2+1)}\ \right).$$

For $\bmu\in \prt^\Phi$, let $L_{\bmu}$, $\bgamma\in \prt_s^{\bmu}$ and $T_{\bgamma}$ be as in Section \ref{LeviTorus}.
Let $\theta$ be a character of $T_{\nu}$.   The \emph{Deligne-Lusztig character} $R_{\nu}(\theta)= R_{T_{\nu}}^{U_n}(\theta)$ is the virtual character of $U_n$ given by
$$\bigl(R_{\nu}(\theta)\bigr)(g)=\sum_{t\in T_\nu\atop \bgamma_t\in \prt_s^{\bmu}} \theta(t) Q_{T_{\bgamma_t}}^{L_{\bmu}}(u),$$
where $g\in c_{\bmu}$ has Jordan decomposition $g=su$ (thus, by Lemma \ref{JordanDecomposition} $C_{U_n}(s)\cong L_{\bmu}$),  and $Q_{T_{\bgamma_t}}^{L_{\bmu}}(u)$ is a Green function for the unitary group (see, for example, \cite{Ca85}).   

Deligne and Lusztig proved
\begin{align*}
C_n&=\{\text{class functions of $U_n$}\}\\
&=\cplx\spanning\{R_\nu(\theta)\ \mid\ \nu\in \prt_n, \theta\in \Hom(T_\nu,\cplx^\times)\},
\end{align*}
so we may define \emph{Deligne-Lusztig induction} by
\begin{equation} \label{DLInduction} \begin{array}{rccl} R_{U_m\oplus U_n}^{U_{m+n}}: & C_m\otimes C_n& \longrightarrow & C_{m+n}\\ &
R_\alpha^{U_m}(\theta_\alpha)\otimes R_\beta^{U_n}(\theta_\beta) & \mapsto & R_{T_\alpha\oplus T_\beta}^{U_{m+n}}(\theta_\alpha\otimes\theta_\beta), \end{array}\end{equation}
for $\alpha\in \prt_m$, $\beta\in \prt_n$, $\theta_\alpha\in \Hom(T_\alpha,\cplx)$, and $\theta_\beta\in \Hom(T_\beta,\cplx)$.

Let $\Lambda$ and $C$ be as in Section \ref{SymmetricFunctions}, except we now give $C$ a graded $\cplx$-algebra structure using Deligne-Lusztig induction.  That is, we define a multiplication $\circ$ on $C$ by  
$$\chi\circ\eta=R_{U_m\oplus U_n}^{U_{m+n}}(\chi\otimes\eta), \quad \text{for $\chi\in C_m$ and $\eta\in C_n$}.$$
We recall the characteristic map defined in Section \ref{EnnolaConj},
$$\begin{array}{rccl} {\rm ch}:& C & \longrightarrow & \Lambda\\ &\pi_{\bmu} & \mapsto & P_{\bmu} \qquad \text{for $\bmu\in \prt^\Phi$}.\end{array}$$
As noted in Proposition \ref{charmap1}, it is immediate that $\ch$ is an isometric isomorphism of vector spaces, but it is not yet clear if $\ch$ is also a ring homomorphism when $C$ has multiplication given by Deligne-Lusztig induction.   

\subsection{The Ennola conjecture II}  \label{CharDLCharacters}  

To prove the Ennola conjecture we require two further ingredients:
\begin{enumerate}
\item[(1)]  Theorem \ref{chDLcharacters} and Corollary \ref{chRingIsomorphism} establish that this new characteristic map is also a ring isomorphism by using a key result by Kawanaka on Green functions to evaluate $\ch(R_\nu(\theta))$.
\item[(2)] Corollary \ref{ennolaconj} uses an explicit decomposition by  Lusztig and Srinivasan of irreducible characters of $U_n$ into Deligne-Lusztig characters to complete the reproof of the Ennola conjecture.
\end{enumerate}

To compute $\ch(R_\nu(\theta))$, we need to write the Green functions $Q_{T_{\bgamma}}^{L_{\bmu}}(u)$ as polynomials in $q$.  These Green functions turn out to be those of the general linear group, except with $q$ replaced by $-q$, which is the essence of Ennola's original idea.  This fact was proven by Hotta and Springer \cite{hotspr} for the case that $p = \mathrm{char}(\fifi)$ is large compared to $n$, and was finally proven in full generality by Kawanaka \cite{kawan2}.

\begin{thm}[Hotta-Springer, Kawanaka] \label{green} The Green functions for the unitary group are given by $Q_{T_{\bgamma}}^{L_{\bmu}}(u) = Q_{\bgamma}^{\bmu} (-q)$, where
$$Q_{\bgamma}^{\bmu}(-q)=\prod_{f\in \Phi_{\bmu}} Q_{\bgamma(f)}^{\bmu(f)}((-q)^{d(f)}),$$
and $Q_\gamma^\mu(q)$ is the classical Green function as in (\ref{PowerToHall}).
\end{thm}

For $\nu=(\nu_1,\nu_2,\ldots, \nu_\ell)\in \prt$ and $\theta=\theta_1\otimes\theta_2\otimes\cdots\otimes \theta_\ell$ a character of $T_\nu$, define
$$p_{\bmu_{\nu\theta}}=\prod_{\varphi\in \Theta} p_{\bmu_{\nu\theta}(\varphi)}(\varphi),\quad\text{where}\quad \bmu_{\nu\theta}(\varphi)=(\nu_i/|\varphi|\ \mid \ \theta_i\in \varphi).$$

\begin{thm} \label{chDLcharacters} Let $\nu=(\nu_1,\nu_2,\ldots,\nu_\ell)\in\prt$, $\theta=\theta_1\otimes\theta_2\otimes\cdots\otimes\theta_\ell$ be a character of $T_\nu$, and $\bnu=\bmu_{\nu\theta}\in \prt^\Theta$.  Then 
$$\ch\bigl(R_\nu(\theta)\bigr)=(-1)^{||\bnu||-\ell(\bnu)}p_{\bnu}.$$
\end{thm}
\begin{proof}  By Theorem \ref{green}, and since $\ch(\pi_{\bmu})=P_{\bmu}$, it suffices to show that the coefficient of $P_{\bmu}$ in the expansion of $(-1)^{|\bnu|-\ell(\bnu)}p_{\bnu}$ is
$$\sum_{t\in T_\nu\atop \gamma_t\in \prt_s^{\bmu}} \theta_\nu(t)Q_{\bgamma_t}^{\bmu}(-q).$$
Since
\begin{equation*}
p_{\bnu}=p_{\nu_1/|\varphi_1|}(\varphi_1)\cdots p_{\nu_\ell/|\varphi_\ell|}(\varphi_\ell),\qquad\text{where} \quad \theta_i\in \varphi_i,
\end{equation*}
the transform (\ref{transform}) implies
\begin{equation*} (-1)^{||\bnu||-\ell(\bnu)}p_{\bnu}=
\prod_{i=1}^\ell \sum_{t_i\in k^\times_{\nu_i}} \theta_i(t_i) p_{\nu_i/d_i}(f_i) = \sum_{t\in T_{\nu}} \theta(t) \biggl(\prod_{i=1}^\ell p_{\nu_i/d_i}(f_i)\biggr),
\end{equation*}
where $f_i\in \Phi$ satisfies $f_i(t_i)=0$ and $d_i=d(f_i)$.   By definition (\ref{gammat}),
\begin{equation}
(-1)^{||\bnu||-\ell(\bnu)}p_{\bnu}=\sum_{t\in T_{\nu}} \theta(t) \biggl(\prod_{f\in \Phi} p_{\bgamma_t(f)}(f)\biggr).
\end{equation}
Change the basis from power-sums to Hall-Littlewood polynomials (\ref{PowerToHall}), and we obtain
\begin{align*}
(-&1)^{||\bnu||-\ell(\bnu)}p_{\bnu}\\ &= \sum_{t\in T_{\nu}} \theta(t) \biggl(\prod_{f\in \Phi} \sum_{\bmu(f)\in \prt_{|\bgamma_t(f)|}} \hspace{-.5cm}Q_{\bgamma_t(f)}^{\bmu(f)}((-q)^{d(f)})(-q)^{-d(f)}{P}_{\bmu(f)}(f ;(-q)^{-d(f)})\biggr)\\
&= \sum_{t\in T_{\nu}} \theta(t)  \sum_{\bmu\in \prt^\Phi\atop \bgamma_t\in\prt_s^{\bmu}} Q_{\bgamma_t}^{\bmu}(-q)P_{\bmu}\\
&=\sum_{\bmu\in \prt^\Phi} \biggl(\sum_{t\in T_\nu\atop \gamma_t\in \prt_s^{\bmu}} \theta(t)Q_{\bgamma_t}^{\bmu}(-q)\biggr)P_{\bmu},
 \end{align*}
as desired.\end{proof}

\begin{cor} \label{chRingIsomorphism} Let multiplication in the character ring $C$ of $U_n$ be given by $\circ$.  Then the characteristic map $\ch: C \rightarrow \Lambda$ is an isometric isomorphism of graded $\cplx$-algebras.
\end{cor}
\begin{proof}  Let $\alpha\in \prt_m$, $\beta\in \prt_n$, $\theta_\alpha\in \Hom(T_\alpha,\cplx^\times)$, $\theta_\beta\in \Hom(T_\beta,\cplx^\times)$, $\balpha=\bmu_{\alpha\theta_\alpha}$, and $\bbeta=\bmu_{\beta\theta_\beta}$.  By (\ref{DLInduction}) and Theorem \ref{chDLcharacters},
\begin{align*} \ch\biggl(R_{U_m\oplus U_n}^{U_{m+n}}\bigl(R_\alpha^{U_m}(\theta_\alpha)\otimes R_\beta^{U_n}(\theta_\beta)\bigr)\biggr)&=\ch\bigl(R_{T_\alpha\oplus T_\beta}^{U_{m+n}}(\theta_\alpha\otimes\theta_\beta)\bigr)\\
 &= (-1)^{m+n-\ell(\balpha)-\ell(\bbeta)}p_{\balpha}p_{\bbeta}\\
&=\ch(R_\alpha^{U_m}(\theta_\alpha))\ch(R_\beta^{U_n}(\theta_\beta)).\end{align*}
Since the Deligne-Lusztig characters span $C$, this product extends linearly to $C$.
\end{proof}

An immediate consequence is that the graded multiplication that Ennola originally defined on $C$ is exactly Deligne-Lusztig induction, or 
\begin{cor} \label{sameprod}  
Let $\chi\in C_m$ and $\eta\in C_n$.  Then
$$\chi\circ \eta=\chi\star\eta.$$
\end{cor}
\noindent We therefore having the advantage of taking either definition when convenience demands.

For $\blam\in \prt^\Theta$, let $L_{\blam}$, $W_{\blam}$, and $T_{\bgamma}$, $\bgamma\in \prt_s^{\blam}$, be as in Section \ref{LeviTorus}.    

Note that the combinatorics of $\bgamma$ almost specifies character $\theta_\gamma$ of $T_{\bgamma}$ in the sense that
$$\theta_{\bgamma}(T_{\bgamma}(\varphi))=\theta_{\varphi}(T_{\bgamma}(\varphi)),\qquad \text{for some $\theta_\varphi\in \varphi$.}$$
In fact, we may define
\begin{equation} 
R_{\bgamma}=R_{T_{\bgamma}}^{U_n}(\theta_{\bgamma})=\ch^{-1}\bigl((-1)^{||\bgamma||-\ell(\bgamma)}p_{\bgamma}\bigr),
\end{equation}
where $\theta_{\bgamma}$ is any choice of the $\theta_{\varphi}$'s.  

For every $\blam\in \prt^\Theta$ there exists a character $\omega^{\blam}$ of $W_{\blam}$ defined by
$$\omega^{\blam}(\bgamma)=\prod_{\varphi\in \Theta} \omega^{\blam(\varphi)}(\bgamma(\varphi)),$$
where $\omega^{\blam}(\bgamma)$ is the value of $\omega^{\blam}$ on the conjugacy class $c_{\bgamma}$ corresponding to $\bgamma\in \prt_s^\Theta$.  

In  \cite{luszsrin}, Lusztig and Srinivasan decomposed the irreducible characters of $U_n$ as linear combinations of Deligne-Lusztig characters, as follows.

\begin{thm}[Lusztig-Srinivasan] \label{LusztigSrinivasan} Let $\blam\in \prt^\Theta_n$.   Then  there exists $\tau^\prime(\blam)\in \intg_{\geq 0}$ such that the class function
$$R(\blam)=(-1)^{\tau^\prime(\blam)+\lfloor n/2\rfloor+\sum_{\varphi\in \Theta} |\blam(\varphi)|+|\varphi|\lceil |\blam(\varphi)|/2\rceil} \sum_{\bgamma\in \prt^{\blam}_s}\frac{\omega^{\blam}(\bgamma)}{z_{\bgamma}}R_{\bgamma}$$
is an irreducible character of $U_n$ ($z_{\bgamma}$ is as in (\ref{zbgamma})).
\end{thm}
\noindent\textbf{Remark.}  The sign
$$(-1)^{\lfloor n/2\rfloor+\sum_{\varphi\in \Theta} |\blam(\varphi)|+|\varphi|\lceil |\blam(\varphi)|/2\rceil}=(-1)^{\fifi\text{-rank of $U_n$ + $\fifi$-rank of $L_{\blam}$}},$$
and Theorem \ref{CharacterDegrees} will show that 
$$\tau^\prime(\blam)=n(\blam^\prime)+||\blam||-\sum_{\varphi\in \Theta} |\varphi| \lceil |\blam(\varphi)|/2\rceil = n(\blam^\prime)+\sum_{\varphi\in \Theta}|\varphi|\lfloor |\blam(\varphi)|/2\rfloor.$$

\begin{cor}[Ennola Conjecture] \label{ennolaconj}
For $\blam\in \prt^\Theta$, there exists $\tau(\blam)\in \intg_{\geq 0}$ such that
$$\left\{\ch^{-1}\biggl((-1)^{\tau(\blam)} s_{\blam}\biggr)\ \mid\ \blam\in \prt_n^\Theta\right\}$$
is the set of irreducible characters of $U_n$.
\end{cor}
\begin{proof}
By Theorem  \ref{LusztigSrinivasan} and Theorem \ref{chDLcharacters},
\begin{align*}\ch&(R(\blam))=(-1)^{\tau^\prime(\blam)+\lfloor n/2\rfloor+\sum_{\varphi\in \Theta} |\blam(\varphi)|+|\varphi|\lceil |\blam(\varphi)|/2\rceil} \sum_{\bgamma\in \prt^{\blam}_s}\frac{\omega^{\blam}(\bgamma)}{z_{\bgamma}}(-1)^{n-\ell(\bgamma)} p_{\bgamma}\\
&=(-1)^{\tau^\prime(\blam)+\lfloor n/2\rfloor+n+\sum_{\varphi\in \Theta}|\varphi|\lceil |\blam(\varphi)|/2\rceil} \sum_{\bgamma\in \prt^{\blam}_s}\frac{\omega^{\blam}(\bgamma)}{z_{\bgamma}}(-1)^{\sum_{\varphi\in \Theta}|\blam(\varphi)|-\ell(\bgamma)} p_{\bgamma}
\end{align*}
Note that the sign character $\omega^{\blam_s}$ of $W_{\blam}$ acts by
$$\omega^{\blam_s}(\bgamma)=(-1)^{\sum_{\varphi\in \Theta}|\bgamma(\varphi)|-\ell(\bgamma)},$$
and that $\omega^{\blam} \otimes \omega^{\blam_s} = \omega^{\blam'}$, so since $\bgamma\in \prt_s^{\blam}$,
\begin{align*}
\ch(R(\blam))&=(-1)^{\tau^\prime(\blam)+\lfloor n/2\rfloor+n+\sum_{\varphi\in \Theta}|\varphi|\lceil |\blam(\varphi)|/2\rceil} \sum_{\bgamma\in \prt^{\blam}_s}\frac{(\omega^{\blam}\otimes \omega^{\blam_s})(\bgamma)}{z_{\bgamma}} p_{\bgamma}\\
&=(-1)^{\tau^\prime(\blam)+\lfloor n/2\rfloor+n+\sum_{\varphi\in \Theta}|\varphi|\lceil |\blam(\varphi)|/2\rceil} \sum_{\bgamma\in \prt^{\blam}_s}\frac{\omega^{\blam^\prime}(\bgamma)}{z_{\bgamma}} p_{\bgamma},
\intertext{and by applying (\ref{SchurToPower}) to a product over $\Theta$,}
&=(-1)^{\tau^\prime(\blam)+\lfloor n/2\rfloor+n+\sum_{\varphi\in \Theta}|\varphi|\lceil |\blam(\varphi)|/2\rceil}s_{\blam^\prime}. & &\qedhere
\end{align*}
\end{proof}

\noindent\textbf{Remark.} There are at least two natural ways to index the irreducible characters of $U_n$ by $\Theta$-partitions: Theorem \ref{LusztigSrinivasan} gives a natural indexing by $\Theta$-partitions, but Corollary \ref{ennolaconj} indicates that the conjugate choice is equally natural.   Since we like to think of Schur functions as irreducible characters, we have chosen the latter indexing.  However, several references, including Ennola \cite{ennola}, Ohmori \cite{ohm}, and Henderson \cite{hender}, make use of the former.

%
%

\section{Characters degrees} \label{SCharacterDegrees}

Let $\blam\in \prt^\Theta$, and suppose $\Box\in \blam$ is in position $(i,j)$ in $\blam(\varphi)$ for some $\varphi\in \Theta$.  The \emph{hook length} $\bh(\Box)$ of $\Box$ is
$$\bh(\Box)=|\varphi|h(\Box),\quad\text{where}\quad h(\Box)=\blam(\varphi)_i-\blam(\varphi)_j^\prime-i-j+1,$$
is the usual hook length for partitions.

\vspace{.25cm}

\noindent\textbf{Example.} For $|\theta|=1$, $|\varphi|=3$ and
$$\blam=\left(\xy<0cm,.3cm>\xymatrix@R=.3cm@C=.3cm{
	*={} & *={} \ar @{-} [l] & *={} \ar @{-} [l]  & *={} \ar @{-} [l] \\
	*={} \ar @{-} [u] &   *={} \ar @{-} [l] \ar@{-} [u] &   *={} \ar @{-} [l] \ar@{-} [u]  &   *={} \ar @{-} [l] \ar@{-} [u] \\
*={} \ar @{-} [u] &   *={} \ar @{-} [l] \ar@{-} [u]  & *={} \ar @{-} [l] \ar@{-} [u]  }\endxy^{(\theta)}, \  \xy<0cm,.5cm>\xymatrix@R=.3cm@C=.3cm{
	*={} & *={} \ar @{-} [l] & *={} \ar @{-} [l] \\
	*={} \ar @{-} [u] &   *={} \ar @{-} [l] \ar@{-} [u]  &  *={} \ar @{-} [l] \ar@{-} [u] \\
	*={} \ar @{-} [u] &   *={} \ar @{-} [l] \ar@{-} [u] \\
	*={} \ar @{-} [u] &   *={} \ar @{-} [l] \ar@{-} [u] }\endxy^{(\varphi)}\ 
\right),\text{ the hook lengths are }\left(\xy<0cm,.3cm>\xymatrix@R=.3cm@C=.3cm{
	*={} & *={} \ar @{-} [l] & *={} \ar @{-} [l]  & *={} \ar @{-} [l] \\
	*={} \ar @{-} [u] &   *={} \ar @{-} [l] \ar@{-} [u] \ar @{} [ul]|{\sss{4}} &   *={} \ar @{-} [l] \ar@{-} [u] \ar @{} [ul]|{\sss{3}}  &   *={} \ar @{-} [l] \ar@{-} [u] \ar @{} [ul]|{\sss{1}}\\
*={} \ar @{-} [u] &   *={} \ar @{-} [l] \ar@{-} [u] \ar @{} [ul]|{\sss{2}} & *={} \ar @{-} [l] \ar@{-} [u]  \ar @{} [ul]|{\sss{1}}}\endxy^{(\theta)}, \  \xy<0cm,.5cm>\xymatrix@R=.3cm@C=.3cm{
	*={} & *={} \ar @{-} [l] & *={} \ar @{-} [l] \\
	*={} \ar @{-} [u] &   *={} \ar @{-} [l] \ar@{-} [u]  \ar @{} [ul]|{\sss{12}} &  *={} \ar @{-} [l] \ar@{-} [u] \ar @{} [ul]|{\sss{3}}\\
	*={} \ar @{-} [u] &   *={} \ar @{-} [l] \ar@{-} [u] \ar @{} [ul]|{\sss{6}}\\
	*={} \ar @{-} [u] &   *={} \ar @{-} [l] \ar@{-} [u] \ar @{} [ul]|{\sss{3}}}\endxy^{(\varphi)}\ 
\right).$$
Example 2 in \cite[I.1]{mac} implies that 
\begin{equation} \label{HookLengthIdentity} \sum_{\Box\in \blam}\bh(\Box)=||\blam||+n(\blam)+n(\blam^\prime).\end{equation}

For $\blam\in \prt^\Theta$, let
$$\eta^{\blam}=\ch^{-1}(s_{\blam}).$$

\begin{thm} \label{CharacterDegrees}
Let $\blam\in \prt^\Theta$ and let $1$ be the identity in $U_{||\blam||}$.  Then
$$\eta^{\blam}(1)=(-1)^{\tau(\blam)} q^{n(\blam^\prime)} \frac{\dd\prod_{1\leq i\leq ||\blam||}(q^i-(-1)^i)}{\dd\prod_{\Box\in \blam}(q^{\bh(\Box)} -(-1)^{\bh(\Box)})},$$
where $\tau(\blam)=||\blam||(||\blam||+3)/2+n(\blam) \equiv \lfloor ||\blam||/2\rfloor+n(\blam) \quad ({\rm mod } \; 2)$.  So for each $\blam$, we have $\chi^{\blam} = (-1)^{\tau(\blam)} \eta^{\blam}$.
\end{thm}
\begin{proof}
We follow the computations in \cite[IV.6]{mac}.  Let $\eta^{\blam}_{\bmu}$ be the value of $\eta^{\blam}$ on the conjugacy class  $c_{\bmu}$.  Then 
$$ s_{\blam} = \sum_{\bmu\in \prt_{||\blam||}^\Phi} \eta^{\blam}_{\bmu} P_{\bmu},$$
implies $\eta^{\blam}_{\bmu} = \langle s_{\blam}, a_{\bmu} P_{\bmu} \rangle$.  Since the $s_{\blam}$ are orthonormal,
\begin{equation} \label{deg}
a_{\bmu} P_{\bmu} = \sum_{\blam} \eta_{\bmu}^{\blam} \bar{s}_{\blam}\qquad \text{where} \qquad 
\bar{s}_{\blam} = \sum_{\bmu\in \prt_{||\blam||}^\Phi} \bar{\eta}^{\blam}_{\bmu} P_{\bmu}\end{equation}
is obtained by taking the complex conjugates of the coefficients of the $P_{\bmu}$.

If $c_{\bmu_1}$ is the conjugacy class corresponding to the identity element $1$ of $U_m$, then $\bmu_1(t-1) = (1^m)$ and $\bmu_1(f) = 0$ for $f \neq t-1$.  From the definition of $a_{\bmu}$, and the fact that $P_{(1^m)}(x;t) = e_m(x)$ (see \cite[III.8]{mac}), we have
\begin{align*} 
a_{\bmu_1} P_{\bmu_1} &= (-1)^m (-q)^{m+m(m-1)}\psi_m(-q^{-1})(-q)^{-m(m-1)/2}P_{(1^m)}(X^{(t-1)};(-q)^{-1}) \\
&= (-1)^m (-q)^{m(m+1)/2} \psi_m((-q)^{-1}) P_{(1^m)} (X^{(t-1)} ; (-q)^{-1})\\
& =  \psi_m(-q) e_m (t-1). \end{align*}
 Therefore, by (\ref{deg}), 
\begin{equation} \label{deg2}
\psi_m(-q) e_m (t-1) = \sum_{||\blam|| = m} \eta^{\blam} (1) \bar{s}_{\blam}.
\end{equation}

Let $\delta:\Lambda \rightarrow \cplx$ be the $\cplx$-algebra homomorphism defined by
$$ \delta(p_m(f)) = \begin{cases} (-1)^{m-1}/(q^m - (-1)^m), & \hbox{if } f = t-1, \\ 0, & \hbox{if } f\neq t-1. \end{cases} $$
Apply $\delta \otimes 1$ to both sides of Equation (\ref{A}) to obtain
$$ \log \Big(\sum_{\blam\in \prt^\Theta} \delta(s_{\blam}) \bar{s}_{\blam} \Big) = \sum_{m\geq 1}\frac{(-1)^{m-1}}{m} p_m(t-1)=\log \prod_i (1 + X_i^{(t-1)}). $$
By exponentiating and expanding the product to $e_m$'s, we have
$$ \sum_{m \geq 0} e_m(t-1) = \sum_{\blam} \delta(s_{\blam}) \bar{s}_{\blam},$$
which gives
$$ e_m(t-1) = \sum_{||\blam|| = m} \delta(s_{\blam}) \bar{s}_{\blam}.$$
Compare coefficients with (\ref{deg2}) to deduce
\begin{equation}\label{delprod}
\eta^{\blam}(1) = \psi_m (-q) \delta(s_{\blam}),
\quad\text{where} \quad
\delta(s_{\blam}) = \prod_{\varphi \in \Theta} \delta(s_{\blam(\varphi)}(\varphi)).
\end{equation}
By the definitions of $\delta$ and $p_m(\varphi)$, $\varphi \in \Theta$, 
$$ \delta(p_m(\varphi)) = (-1)^{|\varphi|m} \big( (-q)^{|\varphi|m} - 1 \big)^{-1} = (-1)^{|\varphi|m} \sum_{i \geq 1} (-q)^{-i|\varphi|m},$$
so
$$\delta(p_m(Y_1^{(\varphi)},Y_2^{(\varphi)},\ldots))=(-1)^{|\varphi|m} p_m((-q)^{-|\varphi|},(-q)^{-2|\varphi|},\ldots).$$

That is, applying $\delta$ to a homogeneous symmetric function $f(Y^{(\varphi)})$ of degree $m$ replaces each $Y_i^{(\varphi)}$ by $(-q)^{-i |\varphi|}$ and multiplies by $(-1)^{|\varphi| m}$.  In particular, for $\lambda \in \prt$, $\varphi \in \Theta$, 
\begin{equation}\label{donSchur} \delta(s_{\lambda}(\varphi)) = (-1)^{|\lambda| |\varphi|} s_{\lambda}\big( (-q)^{-|\varphi|}, (-q)^{-2|\varphi|}, \ldots \big).
\end{equation}
Example 2 in \cite[I.3]{mac} implies
\begin{equation} \label{hook} 
\delta(s_{\lambda}(\varphi)) = (-1)^{|\varphi| |\lambda|} (-q)^{-|\varphi| (|\lambda| + n(\lambda))} \prod_{\Box \in \lambda} \big( 1 - (-q)^{-|\varphi| h(\Box)} \big)^{-1}.
\end{equation}
Combine (\ref{delprod}) and (\ref{hook}) to get
\begin{align*} \delta(s_{\blam}) & = (-1)^{||\blam||}  (-q)^{-||\blam|| - n(\blam)} \prod_{\Box \in \blam} \big( 1 - (-q)^{-\bh(\Box)} \big)^{-1}\\
&=(-1)^{||\blam||} (-q)^{-||\blam|| - n(\blam)+\sum_{\Box\in \blam}\bh(\Box)}\prod_{\Box \in \blam} \big( (-q)^{\bh(\Box)} -1\big)^{-1} \\
&= (-1)^{||\blam||} (-q)^{n(\blam')} \prod_{\Box \in \blam} \big( (-q)^{\bh(\Box)} - 1 \big)^{-1}, \hspace{2cm} \text{(by (\ref{HookLengthIdentity}))}\\
&=  (-1)^{||\blam||+n(\blam^\prime)+\sum_{\Box\in \blam}\bh(\Box)} q^{n(\blam')} \prod_{\Box \in \blam} \big( q^{\bh(\Box)} - (-1)^{\bh(\Box)} \big)^{-1}.
\end{align*}
Since
\begin{equation*}||\blam||+n(\blam^\prime) + \sum_{\Box \in \blam} \bh(\Box)= 
2||\blam||+2n(\blam^\prime)+n(\blam) \equiv n(\blam) \quad (\mathrm{mod}\ 2), 
\end{equation*}
we have
$$\delta(s_{\blam})=(-1)^{n(\blam)}q^{n(\blam^\prime)} \prod_{\Box \in \blam} \big( q^{\bh(\Box)} - (-1)^{\bh(\Box)} \big)^{-1}.$$
Since 
$$\psi_m(-q)=\prod_{i=1}^m (1-(-q)^i)=(-1)^{m+m(m+1)/2}\prod_{i=1}^m(q^i-(-1)^i),$$
we have 
\begin{equation*}
\eta^{\blam}(1)=(-1)^{\tau(\blam)} q^{n(\blam^\prime)}\prod_{1\leq i\leq ||\blam||}(q^i-(-1)^i)\prod_{\Box\in \blam}(q^{\bh(\Box)} -(-1)^{\bh(\Box)})^{-1}.\qedhere\end{equation*}
\end{proof}

By the Littlewood-Richardson rule, for any $\mu, \nu \in \prt$, we have
$$ s_{\mu}s_{\nu} = \sum_{\lambda} c_{\mu \nu}^{\lambda} s_{\lambda}, $$
where $c_{\mu \nu}^{\lambda}$ is the number of tableaux $T$ of shape $\lambda - \mu$ and weight $\nu$ such that the word $w(T)$ is a lattice permutation \cite[I.9]{mac}.  So we have
\begin{equation} \label{schurprod}
s_{\bmu} s_{\bnu} = \sum_{\blam} c_{\bmu \bnu}^{\blam} s_{\blam} \quad \mathrm{where} \quad c_{\bmu \bnu}^{\blam} = \prod_{\varphi \in \Theta} c_{\bmu(\varphi) \bnu(\varphi)}^{\blam(\varphi)}.
\end{equation}

\begin{cor} \label{charprod}
Let $\bmu, \bnu \in \prt^{\Theta}$.  Then $\chi^{\bmu} \circ \chi^{\bnu}$ is a character if and only if every $\blam \in \prt^{\Theta}$ such that $c_{\bmu \bnu}^{\blam} > 0$ satisfies
$$ n(\bmu) + n(\bnu) \equiv n(\blam) + ||\bmu|| \; ||\bnu||\quad (\mathrm{mod}\ 2) $$
\end{cor}
\begin{proof}  Since
$$ \ch(\chi^{\bmu} \circ \chi^{\bnu}) = (-1)^{\tau(\bmu) + \tau(\bnu)} \sum_{\blam} c_{\bmu \bnu}^{\blam} s_{\blam}.$$
the class function $\chi^{\bmu} \circ \chi^{\bnu}$ is a character if and only if every $\blam$ such that $c_{\bmu \bnu}^{\blam} > 0$ satisfies
$$\tau(\bmu) + \tau(\bnu) \equiv \tau(\blam) \quad (\mathrm{mod}\ 2).$$
For each $\blam$ such that $c_{\bmu \bnu}^{\blam} > 0$, we have, for every $\varphi \in \Theta$, that there is a tableau of shape $\blam(\varphi) - \bmu(\varphi)$ and weight $\bnu(\varphi)$.  In particular, we have $|\bmu(\varphi)| + |\bnu(\varphi)| = |\blam(\varphi)|$ for every $\varphi \in \Theta$, so that $||\bmu|| + ||\bnu|| = ||\blam||$.  The result follows from the definition of $\tau$.
\end{proof}

It follows that the operation $\circ$ does not always take a pair of characters to a character.  For example, let $\varphi_1$ be the orbit of the trivial character, and define $\bmu = \bnu = \big( \Box^{(\varphi_1)} \big)$.  Then one $\blam$ such that $c_{\bmu \bnu}^{\blam} > 0$ is 
$$\blam =\Big( \xy<0cm,.3cm>\xymatrix@R=.3cm@C=.3cm{
	*={} & *={} \ar @{-} [l] \\
	*={} \ar @{-} [u] &   *={} \ar @{-} [l] \ar@{-} [u]  \\
*={} \ar @{-} [u] &   *={} \ar @{-} [l] \ar@{-} [u]   }\endxy ^{(\varphi_1)} \Big).$$  
Then we have $n(\bmu) + n(\bnu) = 0$ while $n (\blam) + ||\bmu|| \; ||\bnu|| = 1$, and so $\chi^{\bmu} \circ \chi^{\bnu}$ is not a character by Corollary \ref{charprod}. 

%
%

\subsection{A character degree sum}\label{chardegsums}

Let $d_r$ denote the number of $F$-orbits in $\Theta$ of size $r$.  Since the $F$-action preserves $\mltM_m^*$ and $|\mltM_m^*| = q^m - (-1)^m$, we have
\begin{equation} \label{divsum}
q^m - (-1)^m = \sum_{r|m} rd_r, \qquad \text{for $m\in \intg_{> 0}$}.
\end{equation}

\begin{thm} \label{degreesum}
The sum of the degrees of the complex irreducible characters of $U_m$ is given by
$$\sum_{||\blam|| = m} \chi^{\blam}(1) = (q + 1) q^2 (q^3 + 1) q^4 (q^5 + 1) \cdots\bigl(q^m+\frac{(1-(-1)^m)}{2}\bigr).$$
\end{thm}
\begin{proof}
Following a similar approach to \cite[IV.6, Example 5]{mac}, we consider the coefficient of $t^m$ in the series
$$ S = \sum_{\blam \in \prt^{\Theta}} (-1)^{n(\blam)+||\blam||} \delta(s_{\blam})t^{||\blam||}, $$
where $n(\blam)$ is as in (\ref{sigma}).  Note that $S=S_{\text{o}}S_{\text{e}}$, where  
$$ S_{\text{o}} =\prod_{ |\varphi| \atop \text{odd} } \sum_{\lambda} (-1)^{n(\lambda)} \delta(s_{\lambda}(\varphi)) (-t)^{|\lambda| |\varphi| } \ \text{and}\ S_{\text{e}}=\prod_{ |\varphi|  \atop \text{even} } \sum_{\lambda} \delta(s_{\lambda}(\varphi)) t^{|\lambda| |\varphi| }. $$
Combine the following identity from Example 6 in \cite[I.5]{mac},
$$ \sum_{\lambda} (-1)^{n(\lambda)} s_{\lambda} = \prod_i (1 - x_i)^{-1} \prod_{i < j} (1 + x_i x_j)^{-1} $$
with (\ref{donSchur}) to obtain
\begin{align*} 
\sum_{\lambda \in \prt} (-1)^{n(\lambda)} & (-1)^{|\varphi| |\lambda|} \delta(s_{\lambda}(\varphi)) t^{|\lambda| |\varphi| } \\
& = \prod_{i \geq 1} \big(1 - (t(-q)^{-i})^{|\varphi| } \big)^{-1} \prod_{1 \leq i < j} \big(1 + (t^2(-q)^{-i-j})^{|\varphi| }\big)^{-1}.
\end{align*}
Taking the logarithm, 
\begin{align*} \log S_{\text{o}} &= \sum_{m \text{ odd}} d_m \Big( \sum_{i \geq 1} \sum_{r \geq 1} \frac{(t(-q)^{-i})^{mr}}{r} + \sum_{1\leq i < j} \sum_{r \geq 1} \frac{\big(-(t^2(-q)^{-i-j})^m \big)^r}{r} \Big)\\
&= \sum_{m \text{ odd}} d_m \sum_{r \geq 1}\biggl( \frac{t^{mr}}{r}\sum_{i \geq 1} (-q)^{-imr} + (-1)^r \sum_{1\leq i < j} t^{2mr} (-q)^{-(i+j)mr} \biggr)\\
& = \sum_{m \text{ odd}} d_m\sum_{r \geq 1} \frac{t^{mr}}{r((-q)^{mr} - 1)} \biggl( 1 +(-1)^r \sum_{i \geq 1} t^{mr} (-q)^{-2imr} \biggr)\\
& = \sum_{m \text{ odd}} d_m\sum_{r \geq 1} \frac{t^{mr}}{r(q^{mr} - (-1)^{mr})} \biggl( (-1)^{mr} + \sum_{i \geq 1} t^{mr} q^{-2imr} \biggr),
\end{align*}
where $d_m$ is as in (\ref{divsum}).  Similarly for $|\varphi| $ even, by Example 4 in \cite[I.5]{mac},
$$ \sum_{\lambda} s_{\lambda} = \prod_i (1 - x_i)^{-1} \prod_{i < j} (1 - x_i x_j)^{-1} $$
so
$$ \sum_{\lambda \in \prt} \delta(s_{\lambda}(\varphi)) t^{|\lambda| |\varphi| } = \prod_{i \geq 1} \big(1 - (tq^{-i})^{|\varphi| } \big)^{-1} \prod_{1 \leq i < j} \big(1- (t^2 q^{-i-j})^{|\varphi| } \big)^{-1}.$$
Taking logarithms,
\begin{align*}
 \log S_{\text{e}} &= \sum_{m \text{ even}} d_m \biggl( \sum_{i \geq 1} \sum_{r \geq 1} \frac{(tq^{-i})^{mr}}{r} + \sum_{1 \leq i <j} \sum_{r \geq 1} \frac{(t^2 q^{-i-j})^{mr}}{r} \biggr)\\
&= \sum_{m \text{ even}} d_m\sum_{r \geq 1}  \frac{t^{mr}}{r(q^{mr} - (-1)^{mr})} \biggl((-1)^{mr}+ \sum_{i \geq 1} t^{mr} q^{-2imr} \biggr).
\end{align*}
Set $N = mr$, and by (\ref{divsum}),
\begin{align*} \log S &=\log S_{\text{o}} +\log S_{\text{e}} \\
&= \sum_{N \geq 1} \Big(\frac{(-t)^N}{N} + \sum_{i \geq 1} \frac{t^{2N}q^{-2iN}}{N} \Big)\\
&=\log (1+t)^{-1}+\sum_{i\geq 1} \log(1-t^2q^{-2i})^{-1}
\end{align*}
By exponentiating,
$$ S = \frac{1}{1+t} \prod_{i \geq 1} \frac{1}{1 - t^2 q^{-2i}} = (1 - t) \prod_{i \geq 0} \frac{1}{1 - t^2 q^{-2i}}  = (1 - t) \sum_{l\geq 0} \frac{t^{2l}}{\psi_l(q^{-2})}, $$ 
where the last equality comes from Example 4 in \cite[I.2]{mac}.  Multiply the coefficient of $t^m$ by $(-1)^m|\psi_m(-q)| = (-1)^m\prod_{i = 1}^m (q^i - (-1)^i)$ to obtain
\begin{equation*} \sum_{||\blam|| = m} \chi^{\blam}(1) = (q + 1) q^2 (q^3 + 1) q^4 (q^5 + 1) \cdots\bigl(q^m+\frac{(1-(-1)^m)}{2}\bigr).\qedhere \end{equation*}
\end{proof}

Write $f_{U_m}(q)=\sum_{||\blam||=m} \chi^{\blam}(1)$.  The polynomial $f_{G_m}(q)$ expressing the sum of the degrees of the complex irreducible characters of $G_m = {\rm GL}(m, \fifi)$, was computed in \cite{gow1} for odd $q$ and in \cite{kl} and Example 6 of \cite[IV.6]{mac} for general $q$.  From these results we see that
$$ f_{U_m}(q) = (-1)^{m(m+1)/2} f_{G_m}(-q),$$
another example of Ennola duality.

Gow \cite{gow1} and Klyachko \cite{kl} proved that the sum of the degrees of the complex irreducible characters of $G_n$ is equal to the number of symmetric matrices in $G_n$.  Gow accomplished this by considering the split extension of $G_n$ by the transpose-inverse automorphism,
\begin{equation} \label{gext}
G_n^{+} = \langle G_n, \kappa \mid \kappa^2 = 1, \kappa g \kappa = {^t g}^{-1} \text{ for every } g \in G_n \rangle,
\end{equation}
and showed that every complex irreducible representation of $G_n^+$ could be realized over the real numbers.  Instead, we show directly that the sum of the degrees of the complex irreducible characters of $U_n$ is equal to the number of symmetric matrices in $U_n$, and make conclusions about the reality properties of the characters of $U_n^+$, where
\begin{equation} \label{uext}
U_n^{+} = \langle U_n, \kappa \mid \kappa^2 = 1, \kappa u \kappa = {^t u}^{-1} \text{ for every } u \in U_n \rangle.
\end{equation}
\begin{cor} \label{unsym}
The sum of the degrees of the complex irreducible characters of ${\rm U}(n, \fifisq)$ is equal to the number of symmetric matrices in ${\rm U}(n, \fifisq)$.
\end{cor}
\begin{proof}  Define $S_{U_n}$ to be
$$ S_{U_n} = \{ \text{symmetric matrices in } U_n \}.$$
Then $U_n$ acts on the set $S_{U_n}$ by the action
$$ u \cdot s = {^t u}su, \;\; s \in S_{U_n}, \, u \in U_n $$
To find the number of elements in $S_{U_n}$, we find the stabilizers of the orbits under this action.  Consider the same action of $G_n = {\rm GL}(n, \fifi)$ on its subset of symmetric matrices.  That is, let 
$$ S_{G_n} = \{ \text{symmetric matrices in } G_n \},$$
and let $G_n$ act on the set $S_{G_n}$ by the action 
$$g \cdot s = {^t g}sg, \;\; s \in S_{G_n}, \, g \in G_n.$$  

Gow \cite{gowpre} proved, using the Lang-Steinberg theorem, that there is a one-to-one correspondence between conjugacy classes of $G_n^+$ (\ref{gext}) and $U_n^+$ (\ref{uext}) of elements of the form $g\kappa$ and $u\kappa$, for $g \in G_n$ and $u \in U_n$.  Moreover, Gow proved that this correspondence preserves orders of the elements, and the corresponding centralizers are isomorphic.  The conjugacy classes of order 2 elements of this form correspond to the orbits of symmetric matrices in $G_n$ and $U_n$, and their centralizers to the stabilizers under the action described above.  So to find the stabilizers of the orbits in the case for $U_n$, it is enough to do this for $G_n$.

First let $q$ be odd.  Then the orbits of $S_{G_n}$ correspond exactly to equivalence classes of quadratic forms, of which there are exactly 2 (see, for example, \cite[Chapter 9]{grove} for a complete discussion).  The stabilizers of these orbits are the orthogonal groups corresponding to the quadratic forms.  When $n$ is odd, these two orthogonal groups are isomorphic, with order
$$ |O_1| = |O_2| = 2q^{(n-1)^2/4} \prod_{i = 1}^{(n-1)/2} (q^{2i} - 1).$$
When $n$ is even, the orthogonal groups fixing the two classes of quadratic forms are not isomorphic, and their orders are
$$ |O_1| = 2q^{(n^2 - 2n)/4} (q^{n/2} - 1) \prod_{i = 1}^{(n-2)/2} (q^{2i} -1), \text{ and }$$
$$|O_2| = 2q^{(n^2 - 2n)/4} (q^{n/2} + 1) \prod_{i = 1}^{(n-2)/2} (q^{2i} -1).$$
From Gow's result, the orders of the stabilizers of orbits is the same for $U_n$.  Using the fact that
$$ |U_n| = q^{n(n-1)/2} \prod_{i = 1}^{n} (q^i - (-1)^i),$$
we calculate that the number of elements in $S_{U_n}$ is
\begin{align*}
|S_{U_n}| & = \frac{|U_n|}{|O_1|} + \frac{|U_n|}{|O_2|} \\
    & = (q + 1) q^2 (q^3 + 1) q^4 (q^5 + 1) \cdots\bigl(q^n+\frac{(1-(-1)^n)}{2}\bigr),
\end{align*}
and for this case, the result now follows from Theorem \ref{degreesum}.

When $q$ is even, the orbits of $S_{G_n}$ are not necessarily in one-to-one correspondence with inequivalent quadratic forms.  When $n$ is odd and $q$ is even, however, there is exactly one orbit, and its stabilizer is the unique orthogonal group, $O$ (see \cite[Chapter 14]{grove}), with order
$$|O| = q^{(n-1)^2/4} \prod_{i=1}^{(n-1)/2} (q^{2i} - 1).$$
So the size of $S_{U_n}$ is again found to be
$$|S_{U_n}| = \frac{|U_n|}{|O|} = (q + 1) q^2 (q^3 + 1) q^4 (q^5 + 1) \cdots q^{n-1} (q^n + 1).$$

When $q$ and $n$ are both even, there are two orbits, however these do not correspond to classes of quadratic forms.  One orbit consists of symmetric matrices which have at least one nonzero entry on the diagonal, the order of whose stabilizer, ${\rm Stab}(1)$, is computed in \cite{macwill}.  The order is
$$ |{\rm Stab}(1)| = q^{n^2/4} \prod_{i=1}^{(n-2)/2} (q^{2i} - 1) $$
We note that in \cite{macwill}, the author calls ${\rm Stab}(1)$ the orthogonal group, but this is different from what is commonly defined as the orthogonal group for characteristic 2 in the current literature.

The other orbit, consisting of symmetric matrices with zero diagonal, corresponds to the unique class of alternating form, with stabilizer the symplectic group, $Sp = {\rm Sp}(n, \fifi)$, with order
$$ |Sp| = q^{n^2/4} \prod_{i = 1}^{n/2} (q^{2i} - 1). $$
As in the previous cases, we calculate the size of $S_{U_n}$, which in this case is 
\begin{align*}
|S_{U_n}| & = \frac{|U_n|}{|{\rm Stab}(1)|} + \frac{|U_n|}{|Sp|} \\
    & = (q + 1) q^2 (q^3 + 1) q^4 (q^5 + 1) \cdots (q^{n-1} + 1) q^n. \qedhere
\end{align*} 
 
\end{proof}

The connection between Corollary \ref{unsym} and the characters of $U_n^+$ is most easily understood using twisted Frobenius-Schur indicators, which we now briefly discuss.  More details may be found in \cite{KM} and \cite{bg}.

Let $G$ be a finite group, $\iota$ an order 2 automorphism of $G$, and $(\pi, V)$ a complex irreducible representation of $G$ with character $\chi$.  Suppose that ${^\iota \pi} \cong \hat{\pi}$, where ${^\iota \pi}(g) = \pi({^\iota g})$, and where $\hat{\pi}$ is the contragredient of $\pi$, or the irreducible representation with character $\bar{\chi}$.  This isomorphism gives rise to a nondegenerate bilinear form $B_{\iota}$ satisfying
$$ B_{\iota}(\pi(g) v, \; {^\iota \pi}(g) w) = B_{\iota}(v, w) \hbox{ for all } g \in G \hbox{ and } v, w \in V.$$
This bilinear form is unique up to scalar by Schur's lemma, from which it follows that for every $v, w \in V$, we have
$$ B_{\iota}(v, w) = \epi(\pi) B_{\iota}(w, v), $$
where $\epi(\pi) = \pm 1$.  If ${^\iota \pi} \not\cong \hat{\pi}$, then define $\epi(\pi) = 0$, and we also write $\epi(\chi) = \epi(\pi)$.  Kawanaka and Matsuyama \cite{KM} obtained the following formula for this indicator:
\begin{equation} \label{indic}
\epi(\chi) = \frac{1}{|G|} \sum_{g \in G} \chi(g \; {^\iota g}).
\end{equation}
This formula generalizes that of the classical Frobenius-Schur indicator, which is the case that $\iota$ is trivial, in which case we write $\epi(\pi) = \ep(\pi)$.  Frobenius and Schur proved that $\ep(\pi) = 1$ when $\pi$ is a real representation, $\ep(\pi) = -1$ when $\pi$ is not a real representation but $\chi$ is real-valued, and $\ep(\pi) = 0$ when $\chi$ is not real-valued.  The orthogonality relations and equation (\ref{indic}) imply the following generalization of the Frobenius-Schur involution formula (for a proof, see \cite[Proposition 1]{bg}):
\begin{equation} \label{indsum}
\sum_{\chi \in \text{Irr}(G)} \epi(\chi) \chi(1) = \big| \{ g \in G \mid {^\iota g}^{-1} = g \} \big|.
\end{equation}

For a character $\chi$ of a finite group $H$ which is a subgroup of $G$, we say that $\chi$ {\it extends} to $G$ if there is a character $\chi'$ of $G$ such that $\chi'$ restricted to $H$ is $\chi$.

\begin{cor} \label{extension}
Let $U_n^{+}$ be the split extension of $U_n$ by the transpose-inverse involution, as in (\ref{uext}).  Let $\chi$ be a complex irreducible character of $U_n$.  Then we have the following:
\begin{enumerate}
\item If $\ep(\chi) = 0$, then $\chi$ induces to an irreducible character $\psi = {\rm Ind}_{U_n}^{U_n^+} (\chi)$ of $U_n^{+}$ such that $\ep(\psi) = 1$.

\item If $\ep(\chi) = 1$, then $\chi$ extends to an irreducible character $\chi'$ of $U_n^+$ such that $\ep(\chi') = 1$.

\item If $\ep(\chi) = -1$, then $\chi$ extends to an irreducible character $\chi'$ of $U_n^+$ such that $\ep(\chi') = 0$.
\end{enumerate}
\end{cor}
\begin{proof} If ${^\iota g}={^t g}^{-1}$, then Corollary \ref{unsym} and (\ref{indsum}) imply that $\epi(\chi) = 1$ for every irreducible $\chi$ of $U_n$.  First let $\chi$ of $U_n$ be an irreducible character such that $\ep(\chi) = 0$.  It follows from Frobenius reciprocity that $\psi = {\rm Ind}_{U_n}^{U_n^+}(\chi)$ is irreducible.  Applying the formula for an induced character and (\ref{indic}), we find that $\ep(\psi) = \epi(\chi) = 1$.

Now let $\chi$ be an irreducible character of $U_n$ such that $\ep(\chi) = \pm 1$.  Then $\psi = {\rm Ind}_{U_n}^{U_n^+}(\chi)$ is reducible, and is the sum of two extensions of $\chi$ to $U_n^+$, $\psi = \chi' + \chi''$.  Again, using the formula for an induced character and Equation (\ref{indic}), we calculate that
\begin{align*}
\ep(\chi') + \ep(\chi'') & = \ep(\chi) + \epi(\chi) \\
& = \ep(\chi) + 1.
\end{align*}
If $\ep(\chi) = 1$, then we must have $\ep(\chi') = \ep(\chi'') = 1$.  If $\ep(\chi) = -1$, then $\chi$ cannot have an extension with a Frobenius-Schur indicator equal to 1, and so the only choice is that $\ep(\chi') = \ep(\chi'') = 0$
\end{proof}

\subsection{A second character degree sum}

A $\Theta$-partition $\blam$ is \emph{even} if every part of  $\blam(\varphi)$ is even for every $\varphi\in \Theta$. 
\begin{thm} \label{evendegsum}
The sum of the degrees of the complex irreducible characters of $U_{2m}$ corresponding to $\blam$ such that $\blam'$ is even is given by
$$ \sum_{||\blam|| = 2m \atop \blam' \text{ even} } \chi^{\blam}(1) = (q + 1)q^2(q^3 + 1) \cdots q^{2m-2} (q^{2m-1} +1). $$
\end{thm}
\begin{proof}
To calculate this sum, we find the coefficient of $t^{2m}$ in the series
$$ T = \sum_{\blam' \text{ even}} (-1)^{n(\blam) + ||\blam||} \delta(s_{\blam}) t^{||\blam||} = T_{\text{o}}T_{\text{e}},$$ 
where
$$T_{\text{o}} = \prod_{|\varphi|  \atop \text{odd}} \sum_{\lambda' \text{ even}} (-1)^{n(\lambda)} \delta(s_{\lambda}(\varphi)) t^{|\lambda| |\varphi| } \ \text{and} \ T_{\text{e}} = \prod_{ |\varphi|  \atop \text{even} } \sum_{\lambda' \text{ even}} \delta(s_{\lambda}(\varphi)) t^{|\lambda| |\varphi| }. $$
It follows from the computation in Example 6 of \cite[I.5]{mac} that 
$$ \sum_{\lambda' \text{ even}} (-1)^{n(\lambda)} s_{\lambda} = \prod_{i < j} (1 + x_i x_j)^{-1},$$
and applying this yields
$$\sum_{\lambda' \text{ even}} (-1)^{n(\lambda)} \delta(s_{\lambda}(\varphi)) t^{|\lambda| |\varphi| } = \prod_{1 \leq i < j} \Big( 1 + \big(t^2(-q)^{-i-j} \big)^{|\varphi| } \Big)^{-1}.$$ 
From Example 5(b) in \cite[I.5]{mac}, we have the identity
$$ \sum_{\lambda' \text{ even}} s_{\lambda} = \prod_{i < j} (1 - x_i x_j)^{-1},$$
from which it follows, for $|\varphi| $ even, that
$$ \sum_{\lambda' \text{ even}} \delta(s_{\lambda}(\varphi)) t^{|\lambda| |\varphi| } = \prod_{1 \leq i < j} \Big( 1 - \big(t^2 q^{-i-j} \big)^{|\varphi| } \Big)^{-1}.$$
Taking logarithms gives
\begin{align*} \log T_{\text{o}} &= \sum_{m \text{ odd}} d_m \Big( \sum_{1 \leq i < j} \sum_{r \geq 1} \frac{\big(-(t^2(-q)^{-i-j})^m \big)^r}{r} \Big) \\
& = \sum_{m \text{ odd}} d_m \sum_{r \geq 1} \frac{t^{mr}}{r(q^{mr} - (-1)^{mr})} \Big( \sum_{i \geq 1} t^{mr} q^{-2imr} \Big), \end{align*}
where $d_m$ is as in (\ref{divsum}), and
\begin{align*} \log T_{\text{e}} &= \sum_{m \text{ even}} d_m \Big( \sum_{1 \leq i < j} \sum_{r \geq 1} \frac{\big(t^2 q^{-i-j}\big)^{mr}}{r} \Big)\\
&= \sum_{m \text{ even}} d_m \sum_{r \geq 1} \frac{t^{mr}}{r(q^{mr} - (-1)^{mr})} \Big( \sum_{i \geq 1} t^{mr} q^{-2imr} \Big). \end{align*}
Adding the two terms, making the change $N = mr$, and applying (\ref{divsum}), we obtain
$$ \log T = \sum_{N \geq 1} \sum_{i \geq 1} \frac{t^{2N} q^{-2iN}}{N}, $$
and after exponentiating,
$$ T = \prod_{i \geq 1} \frac{1}{1 - t^2 q^{-2i}} = (1 - t^2) \prod_{i \geq 0} \frac{1}{1 - t^2 q^{-2i}} = (1 - t^2) \sum_{l \geq 0} \frac{t^{2l}}{\psi_l(q^{-2})}.$$
Multiplying the coefficient of $t^{2m}$ by $|\psi_{2m}(-q)| = \prod_{i = 1}^{2m} (q^i - (-1)^i)$, we obtain
\begin{equation*} \sum_{||\blam|| = 2m \atop \blam' \text{ even}} \chi^{\blam}(1) = (q + 1)q^2(q^3 + 1) \cdots q^{2m-2} (q^{2m-1} +1).\qedhere \end{equation*}
\end{proof}

Write $g_{U_m}(q)=\sum_{||\blam||=2m,\blam^\prime\text{ even}}\chi^{\blam}(1)$, and let $g_{G_m}(q)$ denote the corresponding sum for $G_m$.  The polynomial $g_{G_m}(q)$ was calculated in Example 7 of \cite[IV.6]{mac}, and similar to the previous example, we see that we have
$$ g_{U_m}(q) = (-1)^m g_{G_m}(-q). $$

In the case that $q$ is odd, Proposition \ref{evendegsum} follows from the following stronger result obtained by Henderson \cite{hender}.  Let $Sp_{2n} = {\rm Sp}(2n, \fifi)$ be the symplectic group over the finite field $\fifi$.

\begin{thm}[Henderson] \label{symp}
Let $q$ be odd.  The decomposition of ${\rm Ind}_{Sp_{2n}}^{U_{2n}}(\One)$ into irreducibles is given by
$$ {\rm Ind}_{Sp_{2n}}^{U_{2n}}(\One) = \sum_{||\blam|| = 2n \atop \blam' \mathrm{ even}} \chi^{\blam}.$$
\end{thm}

The fact that Proposition \ref{evendegsum} holds for all $q$ suggests that Theorem \ref{symp} should as well.

%
%

\section{A Deligne-Lusztig model}\label{Model}

A {\it model} of a finite group $G$ is a representation $\rho$, which is a direct sum of representations induced from one-dimensional representations of subgroups of $G$, such that every irreducible representation of $G$ appears as a component with multiplicity 1 in the decomposition of $\rho$.  

Klyachko \cite{kl} and Inglis and Saxl \cite{ingsax} obtained a model for ${\rm GL}(n, \fifi)$, where the induced representations can be written as a Harish-Chandra product of Gelfand-Graev characters and the permutation character of the finite symplectic group.  

In this section we show that the same result is true for the finite unitary group, except the Harish-Chandra product is replaced by Deligne-Lusztig induction.  The result is therefore not a model for ${\rm U}(n, \fifi)$ in the finite group character induction sense, but rather from the Deligne-Lusztig point of view.

\subsection{The Gelfand-Graev character}

Let $U_n^\prime={\rm GL}(n, \fldK)^{F^\prime}$ as in (\ref{UnitaryDefinition2}), and let
$$B_{<}=\{u\in U_n^\prime\ \mid\ \text{$u$ unipotent and uppertriangular}\}\subseteq U_n^\prime.$$
Fix a nontrivial character $\psi: \fldK_2^+\rightarrow \cplx^\times$ of the additive group of the field $\fldK_2$ such that the restriction to the subgroup $\{x\in \fldK_2\ \mid\ x^q+x=0\}$ is also nontrivial.   The map
$\psi_{(n)}:B_{<}\rightarrow \cplx$ given by
$$\psi_{(n)}(u)=\psi\bigl(u_{12}+\cdots+u_{\lfloor n/2\rfloor-1,\lfloor n/2\rfloor} + u_{\lfloor n/2\rfloor,\lceil n/2\rceil+1}\bigr),\ \text{for $u=(u_{ij})\in B_<$},$$
is a linear character of $B_<$.   Then 
$$\Gamma_{(n)}^\prime= \Ind_{B_<}^{U_n^\prime}(\psi_{(n)})$$ 
is the \emph{Gelfand-Graev character} of $U_n^\prime$.   Let $\Gamma_{(n)}$ be the corresponding Gelfand-Graev character of $U_n={\rm GL}(n, \fldK)^F$.  For $\blam \in \prt^{\Theta}$, define 
$$\hgt(\blam) = \mathrm{max} \{ \ell(\blam(\varphi)) \ | \ \varphi \in \Theta \}.$$  
The following also appears in Section 5.2 of \cite{ohm}.

\begin{thm} \label{gelgraev}
The decomposition of $\Gamma_{(m)}$ into irreducibles is given by
$$\Gamma_{(m)} = \sum_{\blam\in \prt_m^\Theta \atop \hgt(\blam) = 1} \chi^{\blam}.$$
\end{thm}
\begin{proof}  By \cite[Theorem 10.7]{dellusz},
$$\Gamma_{(m)}=(-1)^{\lfloor m/2\rfloor}\sum_{\bgamma\in \prt_m^\Theta} \frac{(-1)^{m-\ell(\bgamma)}}{z_{\bgamma}} R_{\bgamma}.$$
Therefore, Theorem \ref{chDLcharacters} implies
\begin{align*} 
\ch(\Gamma_{(m)}) & = (-1)^{\lfloor m/2\rfloor}\sum_{\bgamma\in \prt_m^\Theta} \frac{1}{z_{\bgamma}} p_{\bgamma}\\
 			& = (-1)^{\lfloor m/2\rfloor}\sum_{\bgamma\in \prt_m^\Theta} \frac{\omega^{\bgamma_s^\prime}(\bgamma)}{z_{\bgamma}} p_{\bgamma} & & (\text{$\omega^{\bgamma_s^\prime}$ is the trivial character of $W_{\bgamma}$})\\
			&= (-1)^{\lfloor m/2\rfloor}\sum_{\blam\in \prt_m^\Theta\atop \hgt(\blam)=1} s_{\blam} & & (\text{by (\ref{SchurToPower})}).
\end{align*}
Since $\lfloor m/2\rfloor +n(\blam)= \lfloor m/2\rfloor$ for all $\blam$ with $\hgt(\blam)=1$, Theorem \ref{CharacterDegrees} implies
\begin{equation*}\Gamma_{(m)}=\ch^{-1}(\ch(\Gamma_{(m)}))=\sum_{\blam\in \prt_m^\Theta\atop \hgt(\blam)=1} \chi^{\blam}. \qedhere\end{equation*}
 \end{proof}

\subsection{The model}

For a partition $\lambda$, let $o(\lambda)$ denote the number of odd parts of $\lambda$, and for $\blam \in \prt^{\Theta}$, let $o(\blam) = \sum_{\varphi \in \Theta} |\varphi|  o(\blam(\varphi))$. 

\begin{thm} \label{model}
Let $q$ be odd.  For each $r$ such that $0 \leq 2r \leq m$, 
$$\Gamma_{m - 2r} \circ {\rm Ind}_{Sp_{2r}}^{U_{2r}}(\One)  = \sum_{o(\blam') = m - 2r} \chi^{\blam}.$$
Furthermore, 
$$\sum_{0 \leq 2r \leq m} \Gamma_{m - 2r} \circ {\rm Ind}_{Sp_{2r}}^{U_{2r}}(\One) = \sum_{||\blam|| = m} \chi^{\blam}$$
\end{thm}
\begin{proof}
Suppose $\bmu,\bnu\in \prt^\Theta$, such that $\hgt(\bmu) = 1$ and $\bnu'$ is even.   From (\ref{schurprod}), Corollary \ref{ennolaconj}, and Pieri's formula \cite[I.5.16]{mac},
\begin{equation} \label{ModelCharacterProduct} \chi^{\bmu} \circ \chi^{\bnu} = (-1)^{\tau(\bmu) + \tau(\bnu)} \sum_{\blam} \chi^{\blam}, \end{equation}
where the sum is taken over all $\blam$ such that for every $\varphi \in \Theta$, $\blam(\varphi) - \bnu(\varphi)$ is a horizontal $|\bmu(\varphi)|$-strip.  

We now use Corollary \ref{charprod} to show that $\chi^{\bmu} \circ \chi^{\bnu}$ is  a character.  As $\blam(\varphi) - \bnu(\varphi)$ is a horizontal $|\bmu(\varphi)|$-strip, the part $\blam(\varphi)'_i$ is either $\bnu(\varphi)'_i$ or  $\bnu(\varphi)'_i + 1$ for every $i=1,2,\ldots, \ell(\blam(\varphi))$.  
By assumption, $\bnu'$ is even, so $\bnu(\varphi)'_i$ is even for every $\varphi\in \Theta$, and so
$$\binom{\bnu(\varphi)'_i + 1}{2}=\bnu(\varphi)'_i + \binom{\bnu(\varphi)'_i}{2}\equiv\binom{\bnu(\varphi)'_i}{2}\qquad \text{(mod 2).}$$
Thus, $n(\blam(\varphi)) = \sum_i\binom{\blam(\varphi)^\prime_i}{2}\equiv n(\bnu(\varphi)) \; (\mathrm{mod} \ 2).$  The assumption $\hgt(\bmu)=1$ implies $n(\bmu(\varphi)) = 0$, and since $||\bnu||$ is even, 
$$ n(\bmu) + n(\bnu) \equiv n(\blam) + ||\bmu||\,||\bnu|| \quad (\mathrm{mod} \ 2).$$
By Corollary \ref{charprod}, $\chi^{\bmu} \circ \chi^{\bnu}$ is a character.

Use the decompositions of Theorem \ref{symp} and Theorem \ref{gelgraev} in the product (\ref{ModelCharacterProduct}) to observe that the irreducible characters $\chi^{\blam}$ in the decomposition of $\Gamma_{m - 2r} \circ {\rm Ind}_{Sp_{2r}}^{U_{2r}}(\One)$ are indexed by $\blam\in \prt^\Theta_m$ such that for every $\varphi$,  $\blam(\varphi) - \bnu(\varphi)$ is a horizontal $|\bmu(\varphi)|$-strip, where $||\bmu|| = m -2r$, for some $\bnu(\varphi)$ such that $\bnu(\varphi)'$ is even.  Then the number of odd parts of $\blam(\varphi)'$ is exactly $|\bmu(\varphi)|$, and so the $\blam$ in the decomposition must satisfy $\sum_{\varphi \in \Theta} |\varphi|  o(\blam(\varphi)') = ||\bmu|| = m - 2r$.
\end{proof}

\noindent {\bf Acknowledgments.  } The first author would like to thank Stanford University for supplying a stimulating research environment, and, in particular, P. Diaconis for acting as a sounding board for many of these ideas.

The second author thanks the Tata Institute of Fundamental Research for an excellent environment for doing much of this work, and Rod Gow for pointing out that his results in \cite{gowpre} apply to the proof of Corollary \ref{unsym}.

\end{document}